\documentclass[11pt]{amsart}
\usepackage{mathabx} % for '\widecheck' macro
\usepackage{amsmath,amsfonts,amssymb,amsthm,epsfig,color, mathrsfs, scalerel, stackengine} %tikz}
\usepackage{hyperref} % with the backref option, each bibliographical entry indicates where it was cited
\usepackage{bbm}
\usepackage{mathrsfs, euscript}
\usepackage{upgreek}
\usepackage{tensor}
\usepackage{enumitem}

\newcommand{\bb}{\mathbb}

\newcommand{\ZZ}{\bb Z}
\newcommand{\PP}{\bb P}
\newcommand{\R}{\bb R}
\newcommand{\RR}{\bb R}

\newcommand{\NN}{\bb N}

\newcommand{\QQ}{\bb Q}

\newcommand{\hit}{\operatorname{hit}}
\newcommand{\miss}{\operatorname{miss}}
\newcommand{\red}{\operatorname{red}}

\newcommand{\tra}{\operatorname{tr}}
\newcommand{\cl}{\operatorname{cl}}

\newcommand{\ord}{\operatorname{ord}}
\newcommand{\next}{\operatorname{succ}}

\newcommand{\M}{\mathcal M}
\newcommand{\mH}{\mathcal H}

\newcommand{\mZ}{\mathcal Z}
\newcommand{\mT}{\mathcal T}
\newcommand{\U}{\mathcal U}

\newcommand{\Am}{\mathcal A}
\newcommand{\B}{\mathcal B}
\newcommand{\mR}{\mathcal R}
\newcommand{\mV}{\mathcal V}
\newcommand{\mS}{\mathcal S}

\newcommand{\mC}{\mathcal C}

\newcommand{\eL}{\mathcal{L}}
\newcommand{\eJ}{\mathcal{J}}
\newcommand{\mW}{\mathcal{W}}

\newcommand{\mP}{\mathcal{P}}

\newcommand{\diam}{\operatorname{diam}}
\newcommand{\inte}{\operatorname{int}}
\newcommand{\ext}{\operatorname{ext}}

\newcommand{\SL}{\operatorname{SL}}
\newcommand{\GL}{\operatorname{GL}}
\newcommand{\PGL}{\operatorname{PGL}}
\newcommand{\PSL}{\operatorname{PSL}}

\newcommand{\F}{\mathcal F}
\newcommand{\mO}{\mathcal O}
\newcommand{\mN}{\mathcal N}

\newcommand{\Fp}{\mathfrak p}

\newcommand{\Ft}{\mathfrak t}
\newcommand{\Fs}{\mathfrak s}
\newcommand{\Fb}{\mathfrak b}

\newcommand{\sgn}{\operatorname{sgn}}

\newcommand{\Tr}{\operatorname{Tr}}
\newcommand{\fin}{\operatorname{fin}}

\newcommand{\BA}{\operatorname{\bold{Bad}}}

\newtheorem{Theorem}{Theorem}
\numberwithin{Theorem}{section}

\newtheorem{coro}[Theorem]{Corollary}

\newtheorem{theo}[Theorem]{Theorem}
\newtheorem{prop}[Theorem]{Proposition}

\newtheorem{lemm}[Theorem]{Lemma}

\newcommand\reallywidehat[1]{%
\savestack{\tmpbox}{\stretchto{%
  \scaleto{%
    \scalerel*[\widthof{\ensuremath{#1}}]{\kern-.6pt\bigwedge\kern-.6pt}%
    {\rule[-\textheight/2]{1ex}{\textheight}}%WIDTH-LIMITED BIG WEDGE
  }{\textheight}% 
}{0.5ex}}%
\stackon[1pt]{#1}{\tmpbox}%
}

\newtheorem*{lemma*}{Lemma}
\newtheorem*{question*}{Question}

\newtheorem*{theorem*}{Theorem}

\newtheorem*{defn*}{Definition}

\theoremstyle{remark}

\newtheorem{rema}[Theorem]{\sc Remark}
\newtheorem*{rema*}{\sc Remark}

\numberwithin{equation}{section}
\begin{document}
\title[Rational-translation Lagrange equality spectra]{An arithmetic theory of simultaneous rational-translation Lagrange equality spectra}
\author{Jimmy Tseng}

\address{Department of Mathematics and Statistics, University of Exeter, Exeter EX4 4QF, U.K}
\email{j.tseng@exeter.ac.uk}

\begin{abstract}  We show geometric and dimensional properties of simultaneous rational-translation Lagrange equality spectra for both $\BA$ and $\BA(\delta)$ using the {\em projective cocycle steering method,} which we introduce.  In particular, we show the existence of a simultaneous Hall ray for $\BA$, a simultaneous Hall interval for $\BA(\delta)$, regimes of full and positive Hausdorff dimension, a scaling law, empty low spectra in certain cases, a compact set of full Hausdorff dimension in $(2Q\sqrt 3 - \varepsilon, 2Q\sqrt 3)$ where $Q$ is a common denominator of the finite set of rationals giving the translations, and the density of Lagrange equality sets in $\RR$ for Lagrange constants in the Hall ray or Hall interval, respectively.  

Furthermore, we show that the Lagrange equality set for $\BA$, for a finite set of rational translations, is thick but, for the set of all rational translations, is empty.  Likewise, we bound the Hausdorff dimension of the Lagrange equality sets for $\BA(\delta)$ in every nonempty open interval of $\RR$.  For each fixed $Q$, our result on the Hall interval, the density result, and the dimension bounds for $\BA(\delta)$ hold for all sufficiently small $\delta>0$.  Both emptiness results (low spectra and the Lagrange equality set under all rational translations) hold for every $\delta>0$.

%Our results for $\BA(\delta)$ are for all $\delta>0$ sufficiently small depending on $Q$, except our result, which holds for all $\delta>0$, on the emptiness of the Lagrange equality sets under all rational translations.

\end{abstract}

\maketitle

\tableofcontents

\section{Introduction}  Let $x \in \RR \backslash \QQ$, $q \in \NN$, $\delta >0$, and $A \subset \RR$.  Define $\|x\| := \min_{m \in \ZZ} \|x - m\|$.  Let $k(x):= \limsup_{q\rightarrow \infty}(q \|qx\|)^{-1}$ denote the {\em Lagrange constant} or the {\em best constant of Diophantine approximation} of $x$ (\cite[Section~1]{CuFr79},~\cite{EJM25}), $\eL(A):= \{k(x): x \in A\}$, $\BA$ denote the set of badly approximable numbers, and $\BA(\delta)$ denote the set of $\delta$-badly approximable numbers (see~\cite[Definition~1.1]{Ts26}).  The classical Lagrange spectrum is $\eL(\BA)$, which has a filtration given by $\eL\left(\BA(\delta')\right) \subset \eL\left(\BA(\delta)\right)$ for $\delta \leq \delta'$. Now $\BA$ is a winning set and thus, for any countable set $\mathfrak{S} \subset \RR$, we have that $\bigcap_{ \kappa \in \mathfrak{S}} (\BA + \kappa)$ is full Hausdorff dimension (\cite{Sch66}), and, likewise, for any finite set $\mathfrak{S} \subset \RR$, we have that $\bigcap_{ \kappa \in  \mathfrak{S}} (\BA(\delta) + \kappa)$ has positive Hausdorff dimension for all sufficiently small $\delta>0$ (\cite[Theorem~1.13]{Ts26}).  We would, however, like more precise information:  how large are the sets of points whose prescribed translates have the same best constant of Diophantine approximation, and what are the size and structure of the corresponding spectra?  This paper seeks to answer this when $ \mathfrak{S}$ is a finite subset of $\QQ$.  Arithmetic properties and tools (Lemma~\ref{lemm:ApproxCoeffIntegralMob}, strong approximation, the Chinese remainder theorem, the transitivity of subgroup actions on $\PP^1(\ZZ / Q \ZZ)$, etc.) are essential for us.

%In this case, we even have $\BA+r= \BA$ for every $r \in \QQ$.

For $r \in \QQ$, define \[E_r := \{x \in \BA : k(x+r) = k(x)\}, \quad E_r(\delta) := \{x \in \BA(\delta) : k(x+r) = k(x)\} \] and, for a set $F \subset \QQ$, define \begin{align*}
&E_F:= \bigcap_{r \in F} E_r, \quad d_F(t):= \dim_H\left(\eL(E_F) \cap (-\infty,t)\right),  \\  &E_F(\delta):= \bigcap_{r \in F} E_r(\delta), \quad d_{F, \delta}(t):= \dim_H\left(\eL(E_F(\delta)) \cap (-\infty,t)\right) \end{align*}where $\dim_H(\cdot)$ denotes Hausdorff dimension.  We refer to the sets $\eL(E_F)$, $\eL(E_F(\delta))$ as {\em simultaneous rational-translation Lagrange equality spectra} and $E_F, E_F(\delta)$ as their respective {\em Lagrange equality sets}.  Note that $\BA(\delta)$ is well-suited for the Lagrange spectrum because we have that \begin{align}\label{eqn:BadDeltaEndPointBnds}
E_F \cap \{x \in \BA : k(x) < \delta^{-1}\} \subset E_F(\delta) \subset E_F \cap \{x \in \BA : k(x) \leq \delta^{-1}\}. \end{align}  Also note that $E_\ZZ = \BA$ and $E_\ZZ(\delta) = \BA(\delta)$ because $k(x+\ZZ) = k(x)$.  Many properties~\cite[Theorem~1]{Mor18} of the dimension function of the classical Lagrange spectrum $d_{\cl}(t) := \dim_H\left(\eL(\BA) \cap (-\infty,t)\right)$ are retained under our arithmetic theory of simultaneous rational-translation Lagrange equality spectra.  Let  $\Ft_1: = \inf\{t \in \RR: d_{\cl}(t)= 1\}  \approx 3.33$~\cite[Theorem~2.1]{MMPV22} and $s_*:= \frac{1 + 2 \sqrt 2 + \sqrt 5}{2} \approx 3.03$.

\begin{theo}\label{thm:MainBad}  Let $Q \geq 2$ be an integer,  $F:= \frac 1 Q \ZZ$, and $p_Q$ be the largest prime dividing $Q$.  Then we have \begin{enumerate}[label=(\roman*)]
\item\label{it1:HallRay:thm:MainBad} $[4\sqrt 2 \, Q, \infty) \subset \eL(E_F)$,
\item\label{it2:FullHD:thm:MainBad}  $d_F(t) = 1$ for all $t \geq Q \Ft_1$,
\item\label{it3:Scaling:thm:MainBad}  $d_F(t) \geq d_{\cl}(t/Q)$ for all  $t \geq Qs_*$,
\item\label{it4:PosDim:thm:MainBad}  $d_F(t) >0$ for all $t > 3Q$,
\item\label{it5:LowSpec:thm:MainBad}  $\eL(E_F) \cap (-\infty, p_Q) = \emptyset$.

\end{enumerate}
Moreover, \begin{enumerate}[label=(\alph*)]
\item\label{it6:DenE:thm:MainBad}  for every $t\geq4 \sqrt{2} \, Q$, we have $\overline{\{x\in E_F:k(x)=t\}}= \RR$,
\item\label{it7:FullDimCpctSet:MainBad} for every $\varepsilon >0$, there exists a compact set $V_{Q, \varepsilon}$ with $\dim_H V_{Q, \varepsilon}=1$ such that \[V_{Q, \varepsilon} \subset \eL(E_F) \cap (2 Q \sqrt 3 - \varepsilon, 2 Q \sqrt 3),\]
\item for $p$ a prime, we have \[[4\sqrt 2 p, \infty) \subset  \eL\left(E_{\frac 1 p \ZZ}\right) \subset [p, \infty), \quad p \leq \inf  \eL\left(E_{\frac 1 p \ZZ}\right) \leq 3p.\]

\end{enumerate}
\end{theo}
\begin{proof}
 Apply Propositions~\ref{prop:HallRay},~\ref{prop:ScalingLaw},~\ref{prop:PosHD}, and~\ref{prop:FullHDSet} and Corollaries~\ref{cor:FullHD},~\ref{cor:LowSpectra}, and~\ref{cor:QisPStatement}.
\end{proof}
%with $F:=\frac 1 Q \ZZ$

For Lagrange equality sets, we have contrasting sizes for finite $F \subset \QQ$ versus $\QQ$.

\begin{theo}\label{thm:EqualSetsBad}  We have \begin{enumerate}[label=(\roman*)]
\item\label{it1:thickEFinite:thm:EqualSetsBad} for every integer $Q \geq 2$, that  $\dim_H E_{\frac 1 Q \ZZ}=1$ and, moreover, $ E_{\frac 1 Q \ZZ}$ is thick, 
\item\label{it2:emptyEQ:thm:EqualSetsBad} $E_\QQ = \emptyset$.
\end{enumerate}
\end{theo}
\begin{proof}
 Apply Proposition~\ref{prop:EFThick} and Corollary~\ref{cor:EQQIsEmpty}.
\end{proof}

\begin{rema}
Theorems~\ref{thm:MainBad} and~\ref{thm:EqualSetsBad}  have counterparts for $\BA(\delta)$.  See Theorems~\ref{thm:MainBadDelta} and~\ref{thm:EqualSetsBadDelta}, respectively.

\end{rema}

Aside from Theorem~\ref{thm:MainBad}~\ref{it5:LowSpec:thm:MainBad} and Theorem~\ref{thm:EqualSetsBad}~\ref{it2:emptyEQ:thm:EqualSetsBad} and their counterparts for $\BA(\delta)$, every other result in both theorems and their counterparts for $\BA(\delta)$ are proved via the projective cocycle steering method, which we introduce in Section~\ref{sec:SteeringMethod}.  The method, which has a few variants, provides a mechanism for imposing prescribed finite arithmetic conditions on problems that are amenable to symbolic dynamics and may have further applications to this class of problems. Let $R_Q:= \ZZ/Q \ZZ$.  We define a finite, prefix-free alphabet $X$ and its Gauss-Cantor set $K(X)$ in Section~\ref{subsec:ShiftsPress}, its Markov height $\mathfrak{H}(X)$ in Section~\ref{sec:DistWord} , and its associated subgroup $H_Q(X)$ of $\PGL(2, R_Q)$ in the beginning of Section~\ref{sec:SteeringMethod}.   An important special case of one of the variants is transitive projective cocycle steering. 

\begin{theo}\label{prop:transitivesteering}  Let $s \in \RR$, $Q\geq 2$ be an integer, $F \subset \QQ$ be nonempty and finite with $Qr \in \ZZ$ for every $r \in F$, and $X$ be a finite, prefix-free alphabet.  If $\mathfrak{H}(X) < s$ and $H_Q(X)$ acts transitively on $\PP^1(R_Q)$, then $d_F(Qs) \geq \min(1, 2 \dim_H K(X))$.
 \end{theo} 
 
\noindent Theorem~\ref{prop:transitivesteering} is proved in Section~\ref{subsubsec:prooftransitivesteering}.  Finally, we note that, while symbolic dynamics (see the no matching lemma in~\cite{Ts09Non}) can be used for winning sets like $\BA$, the projective cocycle steering method is a stronger tool because it can coordinate Markov heights, Hausdorff dimension, and prescribed finite arithmetic conditions.

\subsection{Preliminaries}
Let $[a_0; a_1, a_2, \cdots]$ denote the unique continued fraction expansion of $x$ and $p_n/q_n$ its $n$-th convergent.  For $p \in \ZZ$ and $q \in \NN$ such that $\gcd(p,q)=1$, we refer to $p/q$ as a {\em reduced rational number} and let \[\Lambda\left(x;p/q\right) := \frac 1 {q^2|x-p/q|}= \frac 1 {q |qx-p|}, \quad  \Lambda_n:=\Lambda(x; p_n/q_n).\]  The quantity $\left(\Lambda\left(x;p/q\right)\right)^{-1}$ is known as the {\em approximation coefficient (of $x$ at $p/q$)}.  We also have (\cite[Section~1]{Mor18}) that \begin{align}\label{eqnPerron} \Lambda_n(x) = [a_{n+1}; a_{n+2}, \cdots] + [ 0;a_n, a_{n-1}, \cdots, a_1]\end{align} and that \begin{align}\label{eqnPerronConv} k(x) = \limsup_{n \rightarrow \infty} \Lambda_n(x).
  \end{align}  Elementary properties of continued fractions can be found, for example, in~\cite{Kh}.  
  
Note that Legendre's theorem (\cite[Theorem~5.4.1]{Reu19}) states if $\Lambda\left(x;p/q\right) >2$, then $p/q$ is a convergent of $x$.  Using Legendre's theorem, the arithmetic structure of our setup allows certain convergents to be mapped to convergents of rational translations.  This correspondence enables us to obtain, in Lemma~\ref{lemm:LagrangeUnderQPush}, a sufficient condition for simultaneous equality of Lagrange constants under rational translations. This arithmetic correspondence is essential to our proofs.

Continued fractions give rise to a (right) cocycle which we will use to construct a group action (see the beginning of Section~\ref{sec:SteeringMethod}).  Let $\mathcal{G}$ be the Gauss map, $\alpha:=[0; a_1, a_2 \cdots]$ where $a_i \in \NN$, and $y \in [0,1]$.  Set $w=(a_1, \cdots, a_n)$, and define \[M_{a_1}:=M_{(\alpha,1)}:=\begin{pmatrix}  0 & 1 \\ 1 & a_1 \end{pmatrix}, \quad M_w:=M_{(\alpha,n)}:=M_{a_1}\cdots M_{a_n}, \quad \mathcal{I}(w):=M_w([0,1]) .\]  Matrices, in this paper, act on numbers (and every element of the Riemann sphere) by their M\"obius actions. By induction on the rule for the formation of convergents and the definition of continued fractions, respectively, we have that \begin{align}\label{eqn:ContFracMatConvergents}
 M_w = \begin{pmatrix}  p_{n-1} & p_n \\ q_{n-1} & q_n \end{pmatrix} \quad \textrm { and } \quad M_w(y) = [0;a_1, \cdots, a_{n-1}, a_n+y] \end{align} where $p_n/q_n = [0;a_1, \cdots, a_n]$.  Note that $M_{(\alpha,n)}$ is a cocycle and it satisfies $M_{(\alpha, n+m)} = M_{(\alpha,n)}M_{(\mathcal{G}^n\alpha,m)}$.  See, for example,~\cite[Section~2]{DK98} and~\cite[Chapter~2]{Via14}.   Note also that, using elementary properties of continued fractions, we have \begin{align}\label{eqn:DerivativeCocycle}
 |M_w'(y)| = 1/(q_n+q_{n-1}y)^2  \end{align}and a bounded distortion estimate \begin{align}\label{eqnCFBndDist}\frac{\sup_{y \in [0,1]} |M_w'(y)|}{\inf_{y \in [0,1]} |M_w'(y)|}= \left(1 + \frac{q_{n-1}}{q_n} \right)^2 \leq 4.
  \end{align}

%For a set $\mS \subset \Sigma$, define $\mS^+ := \{(\theta_0, \theta_1, \cdots): \theta \in \mS\}$.  

We will use shifts extensively.  See, for example,~\cite{LM21} for an introduction.  Let $\Sigma = \NN^\ZZ$ be the set of all bi-infinite sequences of positive integers.  For $\theta:=(\theta_n)_{n \in\ZZ}= \cdots \theta_{-1} \theta_0 \theta_1 \cdots \in \Sigma$, we say that $\theta_n$ is at the {\em $n$-th position} of $\theta$.  Let $\sigma$ denote the shift map $(\sigma\theta)_n = \theta_{n+1}$ and $\mathcal{R}$ denote the reflection map $(\mathcal{R} \theta)_n = \theta_{-n}$.  For integers $\ell < m$, we denote the word from the $\ell$-th position to the $m$-th position of $\theta$ by $\theta |_{[\ell,m]}:=\theta_\ell \cdots \theta_m$.  We also extend this notation to infinite sequences by allowing $\ell = -\infty$ or $m = \infty$.  There is an alternative definition of the Lagrange spectrum due to Perron (see~\cite[Chapter~1]{CF89} or~\cite[Section~1]{Mor18}), namely, if we define \[f(\theta) := [\theta_0; \theta_1, \theta_2, \cdots] + [0;\theta_{-1}, \theta_{-2}, \cdots] \quad \textrm { and } \quad \mathfrak{l}(\theta):= \limsup_{n \rightarrow \infty} f(\sigma^n \theta),\] then $\eL(\BA) = \{\mathfrak{l}(\theta) : \theta \in \Sigma, \mathfrak{l}(\theta) < \infty\}$.  (If $y = [0;\theta_1, \theta_2, \cdots]$, then $k(y + \ZZ) = \mathfrak{l}(\theta)$ for any choice of $\theta_n \in \NN$ for $n \leq 0$.)  Also, let us define \[\mathfrak{m}(\theta) :=\sup_{n\in \ZZ} f(\sigma^n \theta)\] and note that the Markov spectrum is $\{\mathfrak{m}(\theta): \theta \in \Sigma, \mathfrak{m}(\theta) < \infty\}$.  Now, let $n_1< n_2< \cdots$ be a strictly increasing subsequence in $\NN$.  If there exists $\theta(x):=(\theta_j(x))_{j \in \ZZ} \in \Sigma$ such that for every integer $R \geq1$ we have \[a_{n_k+1+j} =\theta_j(x)  \quad \textrm { for } \quad |j|\leq R,\] then we have that \begin{align}\label{eqnPerronSubseqConv} \Lambda_{n_k}(x) \rightarrow f(\theta(x)) \quad \textrm { as } \quad n_k \rightarrow \infty
  \end{align} by (\ref{eqnPerron}).
  
  \subsection{Approximation coefficients under integral M\"obius transformations}  Our first important observation is arithmetic (see also~\cite[Section~3]{LS97}).  When specialized to rational translations, this observation and Legendre's theorem give us our sufficient condition, Lemma~\ref{lemm:LagrangeUnderQPush}.
  
  % (Lemma~\ref{lemm:LagrangeUnderQPush}) to determine numbers with equal Lagrange constants.
  
\begin{lemm}\label{lemm:ApproxCoeffIntegralMob}
Let $x \in \RR \backslash \QQ$ and $M(x):=\frac{ax+b}{cx+d}$ where $a, b, c, d \in \ZZ$ such that $D:=|ad-bc|\geq1$.  Let $p/q\in \QQ$ be a reduced rational number with $cp+dq \neq 0$ and define \[e:= \gcd(ap+bq, cp+dq), \quad \widetilde{P}:= \sgn(cp+dq) \frac{ap+bq}{e}, \quad \widetilde{Q}:=\frac{|cp+dq|}{e}.\]  Then $e | D$ and \begin{align}\label{eqn1:lemm:ApproxCoeffIntegralMob}
 \frac{\Lambda(M(x); \widetilde{P}/\widetilde{Q})}{\Lambda(x;p/q)} = \frac{e^2 q|cx+d|}{D|cp+dq|}. \end{align}  Moreover, along any sequence of reduced rational numbers $p/q \rightarrow x$ such that $cp+dq \neq 0$, we have that \begin{align}\label{eqn2:lemm:ApproxCoeffIntegralMob}
 \frac{\Lambda(M(x); \widetilde{P}/\widetilde{Q})}{\Lambda(x;p/q)} = \frac{e^2}{D}\left(1+o(1)\right). \end{align}
\end{lemm}

\begin{proof} Let $P':= ap+bq$ and $Q':= cp+dq$.  Then we have that $dP' - bQ' = (ad-bc)p$ and $ aQ' - cP' = (ad-bc)q$.   Since $p/q$ is reduced, we have that $e |D$.  Simplifying \[\widetilde{Q}^2\left|M(x) - \widetilde{P}/\widetilde{Q}\right| =\widetilde{Q}^2 \frac{D |qx-p|}{|cx+d| |cp+dq|} = \frac{D  |qx-p||cp+dq|}{e^2  |cx+d| }\] yields (\ref{eqn1:lemm:ApproxCoeffIntegralMob}).  Finally, since $|cp+dq|/q \rightarrow |cx+d|$ as $p/q \rightarrow x$, we have (\ref{eqn2:lemm:ApproxCoeffIntegralMob}).  \end{proof}
  
  %(especially Section~\ref{subsec:ProjCocycleSteer})
  
 \section{Projective cocycle steering method}\label{sec:SteeringMethod} One of our main techniques is the projective cocycle steering method, comprised of projective cocycle steering and two variants, Hall projective cocycle steering and fixed-target projective cocycle steering.  We introduce projective cocycle steering in this section, a transitive specialization in Section~\ref{sec:TransitiveSteering}, the first variant in Section~\ref{sec:HallRay}, and the second variant in Section~\ref{sec:DimEF}.  The method may have further applications.  Let $Q \geq 2$ be an integer for which $Qr \in \ZZ$ for every $r \in F$.  Define $R_Q:=\ZZ / Q \ZZ$ with group of units $R^*_Q$ and \[\PP^1(R_Q) := \{(u,v) \in R^2_Q: \gcd(u,v,Q) = 1\} / R^*_Q\] where $(u,v) \sim (u',v')$ if and only if there exists $\lambda \in R^*_Q$ such that $u=\lambda u'$ and $v = \lambda v'$.  We use $[u:v]$ to denote the equivalence class.  See~\cite[(3.3.1)]{Ste07}. The group $\PGL(2, R_Q):= \GL(2, R_Q)/R^*_Q I$ acts on $\PP^1(R_Q)$ by matrix multiplication.  (Here $I$ denotes the $2 \times 2$ identity matrix.)
 
 Now, if $M=(m_{ij})$ is a matrix with integer entries, then $M \bmod Q = (m_{ij}\bmod Q)$ is a matrix with entries in $R_Q$.  We identify $M$ with its image under this entrywise reduction modulo $Q$.   Consequently, for $w=(a_1, \cdots, a_n)$, we have that $M_w \in \GL(2, R_Q)$ since $\det M_w = \pm 1 \in R^*_Q$.  Define\[\rho_Q(w):= M_w \in \PGL(2, R_Q).\]  We also identify an element of $\GL(2, R_Q)$ with its image in $\PGL(2, R_Q)$ under the canonical projection map.  (Similarly, for $\GL(2, \ZZ)$ and $\PGL(2, \ZZ)$.  Let $m \in \NN$.  For the nonempty, finite set $B=\{w_1, \cdots, w_m\}$ where $w_j = \NN^{r_j}$ for $r_j \in \NN$, $1 \leq j \leq m$, define the subgroup \[H_Q(B) := \langle \rho_Q(w): w \in B\rangle \leq \PGL(2, R_Q).\]

 %\subsection{Approximation coefficients under rational translations}  
 
 \subsection{Matched and unmatched convergents}\label{subsec:MatchedUnmatchedConv}  Recall that the M\"obius action of a $2 \times 2$ upper triangular matrix with $1s$ on the diagonal is translation.  Our next important observation is also arithmetic:  rational translations can come from integral M\"obius transformations.  Consider the family \[A_j := \begin{pmatrix}  1 & j/Q\\ 0 & 1 \end{pmatrix}\begin{pmatrix}  1 & 0\\ Q & Q \end{pmatrix}= \begin{pmatrix}  j+1 & j\\ Q & Q \end{pmatrix},\] which gives the family of rational translations $A_0(x) + j/Q$.  Let $p/q \in \QQ$ be a reduced rational number with $p+q \neq 0$ and let $e_j(p,q)$ denote the $e$ from Lemma~\ref{lemm:ApproxCoeffIntegralMob} for $M(x)=A_j(x)$.  For $a \in \ZZ$ and $b \in \ZZ \backslash \{0\}$, define $\red(a/b)$ to be $a/b$ written in lowest terms with positive denominator.  We have an exact characterization for when $e_j(p,q)$ is maximal.  
 \begin{lemm}\label{lemm:ExactCharMaxTransport}  For a reduced rational number $p/q \in \QQ$, we have that \begin{align}\label{eqn1:lemm:ExactCharMaxTransport}
 e_j(p,q) =Q \iff [p:q]=[-j:j+1] \textrm { in } \PP^1(R_Q). \end{align}  Moreover, let $x \in \RR \backslash \QQ$.  Then along any sequence of reduced rational numbers $p/q \rightarrow x$ such that $[p:q]=[-j:j+1]$ holds, we have that \begin{align}\label{eqn2:lemm:ExactCharMaxTransport}
 \frac{\Lambda\left(A_j(x); \red\left(A_j(p/q) \right)\right)}{\Lambda(x;p/q)} = Q(1 +o(1)). \end{align}

\end{lemm}
 \begin{proof} Note that $e_j(p,q) =Q$ holds if and only if $p+j(p+q) \equiv 0 \pmod Q$ holds.  Thus, if $e_j(p,q) =Q$ holds, we have $p \equiv -j(p+q) \pmod Q$ and $(j+1)(p+q) \equiv -p +p + q \equiv q \pmod Q$.  
 
 We claim $(p+q) \in R^*_Q$.  Assume the claim does not hold.  Let $\ell$ be a prime number such that $\ell | (p+q)$ and $\ell | Q$.  Then, we have that $\ell | p$ and $\ell | q$, which is a contradiction and proves the claim.  Consequently, we have that $[p:q]=[-j:j+1]$.
 
If $[p:q]=[-j:j+1]$ holds, then there exists a $\lambda \in R^*_Q$ such that $p \equiv -j\lambda \pmod Q$ and $q \equiv (j+1)\lambda  \pmod Q$, from which we have that $\lambda \equiv p+q \pmod Q$.  Thus, $p+j(p+q) \equiv 0 \pmod Q$ holds.  This shows (\ref{eqn1:lemm:ExactCharMaxTransport}).  For (\ref{eqn2:lemm:ExactCharMaxTransport}), apply Lemma~\ref{lemm:ApproxCoeffIntegralMob}.  \end{proof}

We say that a convergent $P_n/Q_n$ of $A_j x$ is {\em matched} if there exists a convergent $p_m/q_m$ of $x$ with $p_m+ q_m\neq 0$ such that $P_n/Q_n = \red(A_j(p_m/q_m))$.  Otherwise, it is {\em unmatched}.  Let $\M_j(x):= \{n \in \NN\cup \{0\} : P_n/Q_n \textrm{ is matched}\}$ and $\U_j(x):= \{n \in \NN\cup \{0\} : P_n/Q_n \textrm{ is unmatched}\}$.

\begin{lemm}\label{lemm:UnmatchedConvergents}
Let $x \in \RR \backslash \QQ$ and $j \in \ZZ$.  If $|\U_j(x)|= \infty$, then \[\limsup_{\substack{n \rightarrow \infty \\  n\in \U_j(x)}}\Lambda_n\left (A_j(x)\right)\leq 2Q.\] 
\end{lemm}

\begin{proof}
Assume the conclusion does not hold.  Then there exists infinitely many $n \in \U_j(x)$ and an $\eta >0$ such that $\Lambda(A_j(x); P_n/Q_n) > 2 Q + \eta$.  Recall that \[B_j := \begin{pmatrix}  Q & -j \\ -Q & j+1 \end{pmatrix}\] gives the inverse M\"obius transformation to that of $A_j$.  Let $r_n/s_n = \red(B_j(P_n/Q_n))$.  Since $n \in \U_j(x)$, we have that $r_n/s_n$ is not a convergent of $x$.

Since $B_j$ is continuous away from its pole (at $(j+1)/Q$), we have that $r_n/s_n$ is a sequence of reduced rational numbers for which Lemma~\ref{lemm:ApproxCoeffIntegralMob} applies.  This yields \[\frac{\Lambda(x;r_n/s_n )}{\Lambda(A_j(x);P_n/Q_n )} \geq \frac 1 Q (1 - o(1)),\] which, for $n$ large enough, Legendre's theorem gives a contradiction. 

%Since $P_n/Q_n \rightarrow A_j(x)$, we have \[\frac{|-Q P_n +(j+1) Q_n|}{Q_n} \rightarrow |-Q A_j(x) + j+1| =\frac 1 {|x+1|}\] as $n \rightarrow \infty$.
% 
% 
% as $n \rightarrow \infty$
%Let $e_k := \gcd (Q P_n - j Q_n, -QP_n+(j+1) Q_n) = \gcd(QP_n, Q_n) = \gcd(Q, Q_n)$ since $P_n/Q_n$ is a reduced rational number.    
\end{proof}
 
 \begin{lemm}\label{lemm:LagrangeUnderQPush}  Let $\lambda \geq \sqrt{5}$, $J \subset R_Q$, and $x \in \RR \backslash \QQ$ such that $k(x) = \lambda$.  Let $p_n/q_n$ denote the $n$-th convergent of $x$.  If, for every $j \in J$, there exists an infinite set $N_j \subset \NN$ such that, for all $n \in N_j$, we have \[[p_n:q_n] = [-j:j+1] \textrm { in } \PP^1(R_Q) \quad \textrm { and } \quad \lim_{\substack{n \rightarrow \infty \\ n \in N_j}}\Lambda_n(x) =\lambda,\] then $k(A_j(x))=Q\lambda$ for every $j \in J$.
 
\end{lemm}
\begin{rema*}
 Since the Lagrange constant is invariant under integral translations, $k(A_j(x))$ is well-defined for $j \in R_Q$.   

\end{rema*}
 \begin{proof} Fix $j \in J$.  By Lemma~\ref{lemm:ExactCharMaxTransport}, we have, as $n \rightarrow \infty$ along $N_j$, \[\Lambda(A_j(x), \red(A_j(p_n/q_n))) = Q(1 +o(1))\Lambda_n(x) \rightarrow Q \lambda >2.\]  Legendre's theorem implies that, for all $n \in N_j$ large enough, $ \red(A_j(p_n/q_n))$ is a convergent of $A_j(x)$.  Consequently, $k(A_j(x)) \geq Q \lambda$ by (\ref{eqnPerronConv}).
 
 For the reverse inequality, first consider the matched convergents.  Any infinite sequence of matched convergents $P_m/Q_m$ where $m \in \M_j$ comes from an infinite sequence of convergents $p_{n(m)}/q_{n(m)}$ of $x$.  Lemma~\ref{lemm:ApproxCoeffIntegralMob} gives that \[\Lambda_{m}(A_j(x)) = \frac{e^2_j(p_{n(m)},q_{n(m)})}{Q}(1 +o(1)) \Lambda_{n(m)}(x) \leq Q (1 +o(1)) \Lambda_{n(m)}(x).\]  Hence, by (\ref{eqnPerronConv}), we have \[\limsup_{\substack{m \rightarrow \infty \\  m\in \M_j(x)}}\Lambda_{m}(A_j(x))\leq Q\lambda.\]  Finally, for the unmatched convergents, apply Lemma~\ref{lemm:UnmatchedConvergents} if $|\U_j(x)| =\infty$.  If $|\U_j(x)| < \infty$, the unmatched convergents do not contribute.  This yields the desired result. \end{proof}

\subsection{Complete shifts associated with finite prefix-free alphabets and pressure}\label{subsec:ShiftsPress}  We return to considering shifts.  For $m \in \NN$, the nonempty, finite set $B=\{w_1, \cdots, w_m\}$, where $w_j \in \NN^{r_j}$ for $r_j \in \NN$, $1 \leq j \leq m$ and $w_i \neq w_j$ for $ i \neq j$, is called a {\em finite alphabet}. The $w_j$ are words over $\NN$ and also letters of the alphabet $B$.  Note that we reserve the term {\em digit} for a letter of the alphabet $\NN$ to emphasize that a digit is a natural number.  If all the digits of all the $w_j$ are further restricted to be in $\{1,2\}$, then we refer to $B$ as a {\em finite binary alphabet}. For $w_j =: (w_{j,1}, \cdots, w_{j, r_j})$, we say that $w_{j,r}$ is at the {\em $r$-th position} of $w_j$.  Let $|\cdot|_B$ denote the  {\em length over $B$ of a word over $B$} and $| \cdot|_{\NN}$ denote the {\em length over $\NN$ of a word over $\NN$}.  Thus, for $n \in \NN$ and $w:=(\beta_1, \cdots, \beta_n)$ where $\beta_j \in B$ for $1 \leq j \leq n$, we have that $|w|_B=n$ and $|w|_{\NN} = \sum_{j=1}^n |\beta_j|_{\NN}$.  For concatenations, we follow the convention in~\cite[Section~3]{Mor18} and let juxtaposition also denote concatenation and powers to denote repeated concatenation:  if $\cdots, w_{-1}, w_0, w_1, \cdots$ are words over $\NN$, then $\cdots w_{-1} w_0 w_1 \cdots$ is their concatenation and, for $N \in \NN$, $w_0^N$ is the $N$-fold concatenation of $w_0$.  More generally, we denote the set of all concatenations of exactly $N$ words of $B$ by \[B^N:= \{\beta_1 \cdots \beta_N: \beta_i \in B\}\] and let $B^+ := \bigcup_{N=1}^\infty B^N$ and $B^*$ be the union of $B^+$ with the empty word.

We now associate to $B$ the {\em Gauss-Cantor set} $K(B) \subset [0,1]$ defined by Moreira in~\cite[Section~3]{Mor18}:  \[K(B):= \{[0;\beta_1, \beta_2, \cdots]: \beta_i \in B \textrm { for all } i\geq 1\}.\]  Here, we also allow the degenerate case:  if $B=\{b\}$, then $K(B)=\{[0;b,b, \cdots]\}$.   For a nonempty word $w = (a_1,\cdots, a_n)$ over $\NN$, define \[w^{\mathsf T} := (a_n, \cdots, a_1) \quad \textrm { and } \quad B^{\mathsf T} :=\{w^{\mathsf T}: w \in B\}.\]  Similarly for $\theta:=(\theta_1, \theta_2, \cdots) \in \NN^\NN$, define $\theta^\mathsf{T}:= (\cdots, \theta_2, \theta_1)$.  Note that we have $K(B) = K(B^N)$ and, by \cite[Proposition~2]{Mor18}, that \begin{align}\label{eqnHDKHDKTrans} \dim_H K(B) = \dim_H K(B^{\mathsf T}). \end{align}  Now, for nonempty $Z \subset \Sigma$, define \[\mathcal{K}(Z) := \{[0;\theta_0, \theta_1, \cdots] : \theta \in Z\}.\]

A word $w$ is a {\em prefix} of a word or one-sided infinite sequence $v$ if $v=wu$ for some word or one-sided infinite sequence $u$, respectively.  Furthermore, if $w_i$ does not begin with $w_j$ for $i\neq j$, then the finite alphabet $B$ is called {\em prefix-free}.   Let $A$ be a finite, prefix-free alphabet.  Let $\Sigma^0(A)$ be the set of bi-infinite sequences $\theta$ obtained by bi-infinite concatenations of the letters of $A$ such that the $0$-th position of $\theta$ is the first position of some letter of $A$.  The {\em complete shift associated with $A$} is the compact subset \[\Sigma(A) := \{\sigma^s \theta: \theta \in \Sigma^0(A), s \in \ZZ \}\] of $\Sigma$.  Note that we have $\Sigma(A) = \Sigma(A^N)$.  Since $A$ is prefix-free, the word $w:=\alpha_1 \cdots \alpha_n$ has a unique factorization into the elements of $A$, and thus we have $|w|_A =n$.  Let us also define $\Sigma^+(A)$ to be the set of infinite sequences obtained by infinite concatenations of the letters of $A$.  Finally, note that we also have an total ordering on $A$ as follows:  for $\alpha, \alpha' \in A$, we say that $\alpha < \alpha'$ if and only if $M_\alpha(1/2) < M_{\alpha'}(1/2)$.  Since $A$ is prefix-free, we have that $\inte(M_\alpha([0,1])) \cap \inte(M_{\alpha'}([0,1])) = \emptyset$ whenever $\alpha \neq \alpha'$, and thus $<$ is a strict total ordering on $A$.

\begin{lemm}\label{lemm:AlphabetShiftSameCantorDim} We have $\dim_H \mathcal{K}(\Sigma(A)) = \dim_H \mathcal{K} (\mathcal {R}\Sigma(A))= \dim_H K(A).$
\end{lemm}
\begin{proof}
 Let $\theta \in \Sigma(A)$.  Then there exists $w:=a_1 \cdots a_{|w|_{\NN}} \in A$ such that the $0$th-position of $\theta$ is the $(j+1)$-st position of $w$ for some $j = 0, \cdots, |w|_{\NN}-1$.  Consequently, \[\mathcal{K}(\Sigma(A)) = \bigcup_{w \in A} \bigcup_{j=0}^{|w|_{\NN}-1}M_{a_{j+1} \cdots a_{|w|_{\NN}}}(K(A)), \quad \mathcal{K} (\mathcal {R}\Sigma(A))=\bigcup_{w \in A} \bigcup_{j=0}^{|w|_{\NN}-1}M_{a_{j+1} \cdots a_1}(K(A^\mathsf{T})).\]  Since the unions are finite, the mappings are bilipschitz on $[0,1]$, and (\ref{eqnHDKHDKTrans}) holds, the desired results follow.
\end{proof}

The set $K(A)$ associated to a finite, prefix-free alphabet $A$ when $|A|\geq2$ is a $C^2$-regular Cantor set and is non-essentially affine and the same applies to its image under a $C^2$-diffeomorphism (for which $K(A)$ is contained in an open set in the domain of the diffeomorphism)~\cite[Section~3 and Proposition~1]{Mor18}.  To such sets, we will apply 

\begin{lemm}[Moreira's dimension formula, {\cite[Theorem, Section~2]{Mor18}}]\label{lemm:Moreira}  If $K$ and $K'$ are $C^2$-regular Cantor sets and $K$ is non-essentially affine, then \[\dim_H(K + K') = \min(\dim_H K+ \dim_H K', 1).\]
 
\end{lemm}

For the finite, prefix-free alphabet $A$, we have (\cite[Section~1]{MaUr99}) that the cocycles $\{M_w: w \in A\}$ (when $ w \neq 1$) form a regular conformal iterated function system (IFS) on $[0,1]$.  (Note that prefix-free gives the open set condition.)  Let $u \geq 0$ and $N \in \NN$.  Define \begin{align}\label{eqn:PartitionSumsPressure}
 Z_N^-(A,u) := \sum_{w \in A^N}\inf_{y \in [0,1]} |M_w'(y)|^u, \quad Z_N^+(A,u) := \sum_{w \in A^N}\sup_{y \in [0,1]} |M_w'(y)|^u, \end{align} and, using (\ref{eqnCFBndDist}), the pressure \[P_A(u) := \lim_{N \rightarrow \infty} \frac 1 N \log Z_N^+(A, u) = \lim_{N \rightarrow \infty} \frac 1 N \log Z_N^-(A, u).\]  Note that $P_{A^m}(u) = m P_A(u)$ for any $m \in \NN$.  We have a standard consequence of pressure (see~\cite{MaUr03} for an introduction).

\begin{lemm}\label{lemm:pressure}  Let $A$ be a finite, prefix-free alphabet and $u \geq 0$.  \begin{enumerate}[label=(\roman*)]
\item\label{it1:lemm:pressure} If $u < \dim_H K(A)$, then $P_A(u) >0$ and, there exists an $\eta>0$ and $N_0 \in \NN$ (both depending on $u$), such that, for all $N \geq N_0$, we have that $Z_N^-(A, u)\geq e^{\eta N}$.
\item\label{it2:lemm:pressure} If \[\sum_{w \in A} \inf_{y \in [0,1]} |M_w'(y)|^u>1,\] then $\dim_H K(A) > u$.

\end{enumerate}

\end{lemm}
\begin{proof}  First, we assume that the word $1 \notin A$.   We prove \ref{it1:lemm:pressure}.  Let $A^\NN := \{b_1b_2\cdots: b_i \in A\}$ and $\pi$ be the map in~\cite[Section~1]{MaUr99}). Then we have \[K(A) = \pi(A^\NN).\]  Since our cocycles form a regular conformal IFS, \cite[Theorem~1.1]{MaUr99} gives the pressure conclusion, from which immediately follows the conclusion for the exponential bound.   

Next, we prove \ref{it2:lemm:pressure}.  Let $w = w_1 \cdots w_N$ where $w_n \in A$ for $1 \leq n \leq N$.  By the (repeated use of) the chain rule, we have that \[\inf_{y \in [0,1]} |M_w'(y)|^u \geq \prod_{n=1}^N \inf_{y \in [0,1]} |M_{w_n}'(y)|^u,\] from which it follows that \[Z_N^-(A,u) \geq  \sum_{(w_1, \cdots, w_N) \in A^N}\prod_{n=1}^N \inf_{y \in [0,1]} |M_{w_n}'(y)|^u = \left(\sum_{w \in A}\inf_{y \in [0,1]}|M_{w}'(y)|^u  \right)^N.\]  Note the equality follows by repeated distribution.  Hence, the result follows by~\cite[Theorem~1.1]{MaUr99}.

Finally, we assume that the word $1 \in A$.  Then we replace $A$ with $A^2$.  Since $K(A) = K(A^2)$ and $P_{A^2}(u) = 2 P_A(u)$, the pressure conclusion of \ref{it1:lemm:pressure} follows.  For the exponential bound, the above proof for $A^2$ gives the exponential bound on for $Z_{2N}^-(A, u)$.  The chain rule argument in \ref{it2:lemm:pressure} generalizes to show that $Z_{2N+1}^-(A,u) \geq Z_{2N}^-(A, u) Z_1^-(A,u)$, giving the exponential bound conclusion of \ref{it1:lemm:pressure} for the case $1 \in A$.  For \ref{it2:lemm:pressure}, the chain rule shows that if hypothesis of \ref{it2:lemm:pressure} holds for $A$, then it also holds for $A^2$.  Since $K(A) = K(A^2)$ holds, the conclusion of \ref{it2:lemm:pressure} holds for the case $1 \in A$.
\end{proof}

Let $B$ be a finite alphabet (not necessarily prefix-free), $\widetilde{H}$ be finite group, and $\rho: B \rightarrow \widetilde{H}$ be a mapping.  Extend $\rho$ to any nonempty word $w = (w_1, \cdots, w_N)$ where $w_n \in B$ for $1 \leq n \leq N$ by defining $\rho(w) = \rho(w_1) \cdots \rho(w_N)$.  Define the subgroup $H := H_{B, \rho}:=\langle \rho(w) : w \in B \rangle$.  Since $H$ is finite, we have that, for any $h \in H$, there exists a $w \in B^+$ such that $h= \rho(w)$.  Note that $1_H = \rho(w^{\ord(\rho(w))})$ for $w \in B^+$ where $\ord(\rho(w))$ denotes the order of $\rho(w)$ in $H$.  We also use the convention that $\rho$ applied to the empty word is equal to $1_H$.

\begin{lemm}\label{lemm:IdentityAlphabet}  Let $Y$ be a finite, prefix-free alphabet such that $\dim_H K(Y) >0$ and $H:= H_{Y, \rho}$.  Then, for every $\varepsilon>0$, there exists a finite alphabet $A$ consisting of concatenations of the elements of $Y$ such that \begin{enumerate}[label=(\roman*)]
\item\label{eqn1:lemm:IdentityAlphabet}  $|a|_Y=|a'|_Y$ and $|a|_{\NN} = |a'|_{\NN}$ for every $a, a' \in A$, 
\item $\rho(a)= 1_H$ for every $a \in A$,
\item $\dim_H K(A) > \dim_H K(Y) - \varepsilon.$
\end{enumerate}
 
\end{lemm}
\begin{rema*}
Note that the second equality in \ref{eqn1:lemm:IdentityAlphabet} implies $A$ and $A^{\mathsf T}$ are prefix-free.
\end{rema*}

\begin{proof}  Let $d = \dim_H K(Y)$.  Choose $\max(0, d - \varepsilon) < u < d$.  By Lemma~\ref{lemm:pressure}, there exists an $\eta>0$ and $N_0 \in \NN$ (both depending on $u$), such that, for all $N \geq N_0$, we have that $Z_N^-(Y, u)\geq e^{\eta N}$.  Let \[L_{\min}:=\min_{\gamma \in Y}|\gamma|_\NN, \quad L_{\max} :=  \max_{\gamma \in Y}|\gamma|_\NN, \quad  \Delta:= L_{\max} - L_{\min}, \quad \varphi_u(w) :=\inf_{y \in [0,1]}|M'_w(y)|^u.\]  Partition $Y^N$ using the equivalence relation $w \sim w'$ if and only if $\rho(w) =\rho(w')$ and $|w|_\NN = |w'|_\NN$.  There are at most $|H|(N\Delta+1)$ equivalence classes.  Now for some pair $(h_N, \ell_N) \in H \times \{NL_{\min}, NL_{\min}+1, \cdots, N L_{\max}\}$, the equivalent class \[\mC_N := \{w \in Y^N: \rho(w) = h_N, |w|_\NN = \ell_N\}\] satisfies \begin{align}\label{eqn2:lemm:IdentityAlphabet} \sum_{w \in \mC_N}  \varphi_u(w) \geq \frac {Z_N^-(Y, u)}{|H|(N\Delta+1)}.
  \end{align}
  
  For each $h \in H$, there exists a $c_h \in Y^+$ such that $\rho(c_h) = h^{-1}$.  Let \[\kappa_u := \min_{h \in H} \varphi_u(c_h) >0, \quad A_N := \{w c_{h_N}: w \in \mC_N \}.\]  Here the strict inequality follows by (\ref{eqn:DerivativeCocycle}) and the elements of $A_N$ are distinct.  Hence, we have that $\rho(w c_{h_N}) = h_N h^{-1}_N = 1_H$, $|w c_{h_N}|_Y = N+ |c_{h_N}|_Y$, and $|w c_{h_N}|_\NN = \ell_N+ |c_{h_N}|_\NN$ for all $w \in \mC_N$.  By the chain rule and (\ref{eqn2:lemm:IdentityAlphabet}), we have that \begin{align*}\sum_{w \in A_N} \varphi_u(w) \geq \frac {\kappa_u Z_N^-(Y, u)}{|H|(N\Delta+1)} \geq \frac {\kappa_u e^{\eta N}}{|H|(N\Delta+1)} .
  \end{align*}  (Note that, since $c_{h_N}$ is a fixed word for $\mC_N$, summing over $A_N$ is equivalent to summing over $\mC_N$.)  For large enough $N$, we have $\dim_H(A_N) > u > d- \varepsilon$ by Lemma~\ref{lemm:pressure}.  Taking $A= A_N$ for such an $N$ yields the desired result.\end{proof}

\begin{lemm}\label{lemm:HighXRemoveExtremes}  Let $X$ be a finite prefix-free alphabet such that $\dim_H K(X) >0$ and $H:= H_{X, \rho}$. Then, for every  $R \geq 1$, $\varepsilon>0$, and $N \in \NN$, there exists an integer $n \geq N$ such that, for every $E \subset X^n$ with $|E|\leq R$, we have that \[\langle \rho(w) : w \in X^n \backslash E \rangle = H \quad \textrm{ and } \quad \dim_H K(X^n \backslash E) > \dim_H K(X) - \varepsilon.\]
 
\end{lemm}
\begin{proof}  Choose $L \in \NN$ such that $|X|^L > |H|$.  Then there exist distinct $u, v \in X^L$ such that $\rho(u) = \rho(v)=h$.  Let $U := u^{\ord{h}}$ and $V:= v^{\ord{h}}$.  Then $U \neq V$ but $|U|_X = |V|_X= \ord(h) L$ and $\rho(U) = \rho(V)= 1_H$ both hold.  Thus, for any integer $k \geq 1$ and $w \in X$, the set \[\mS := \{W_1 \cdots W_k w : W_i \in \{U, V\}\}\] consists of $2^k$ distinct words such that, for any $\gamma, \gamma' \in \mS$ we have that $|\gamma|_X = |\gamma'|_X = \ord(h) kL+1$ and $\rho(\gamma) = \rho(\gamma') = \rho(w)$.

Let $d:=\dim_H K(X)$.   Now choose $\max(0, d - \varepsilon) < s < d$.  By Lemma~\ref{lemm:pressure}, there exists $n_0 \in \NN$ such that $Z_n^-(X, s ) > R+1$ for every $n \geq n_0$.  Choose $k$ so that $2^k \geq R+1$ and $\ord(h) kL+1\geq \max(n_0,N)$ both hold.  Set $n:= \ord(h) kL+1$.  Consequently, for each $w \in X$, there exists at least $R+1$ distinct words $\gamma$ in $X^n$ such that $\rho(w) = \rho(\gamma)$.  Since $|E| \leq R$, then the first conclusion holds.

Now, by (\ref{eqn:DerivativeCocycle}), we have that \[\sum_{\gamma \in E} \inf_{y \in [0,1]}|M'_\gamma(y)|^s\leq |E|\leq R,\] which gives that \[\sum_{\gamma \in X^n \backslash E} \inf_{y \in [0,1]}|M'_\gamma(y)|^s \geq Z_n^-(X, s) - R> 1.\] Lemma~\ref{lemm:pressure} yields the second conclusion and concludes the proof.\end{proof}

\subsection{Distinguished-word coding}\label{sec:DistWord}  In this section, we apply Moreira's dimension formula to our setup.  The geometric core of our distinguished word coding is inspired by~\cite[Lemma~3]{Mor18}.  Let $Y$ be a finite prefix-free alphabet  and \[\mathfrak{H}(Y) := \max_{\theta \in \Sigma(Y)}f(\theta),\] which is well-defined because $f$ is continuous on the compact set $\Sigma(Y)$.  We refer to $\mathfrak{H}(Y)$ as the {\em (Markov) height (of $Y$)}.  Let $A$ be a finite prefix-free alphabet consisting of concatenations of the elements of $Y$ such that $A^{\mathsf T}$  is also prefix-free.  Let $\widetilde{n}_-, \widetilde{n}_+ \in \NN \cup \{0\}$, define a fixed word $\widetilde{W} := a_{-\widetilde{n}_-} \cdots a_0 \cdots a_{\widetilde{n}_+}$ where $a_j \in \NN$ for $-\widetilde{n}_- \leq j \leq \widetilde{n}_+$., and refer to $a_0$ as the {\em distinguished digit (of $\widetilde{W}$)}.  For $\Gamma := \gamma_1 \gamma_2 \gamma_3 \cdots, \Delta:= \delta_1 \delta_2 \delta_3 \cdots \in \Sigma^+(A)$, define the element of $\Sigma$ \[\Theta_{\widetilde{W}}(\Gamma, \Delta) := \cdots \delta_3 \delta_2 \delta_1 \widetilde{W} \gamma_1 \gamma_2 \gamma_3 \cdots\] whose $0$-th position is the distinguished digit $a_0$ and let $\mS_{\widetilde{W}} := \{\Theta_{\widetilde{W}}(\Gamma, \Delta) : \Gamma, \Delta \in \Sigma^+(A)\}.$  We say $\mS_{\widetilde{W}}$ is {\em uniformly separated above $\Sigma(Y)$} if there exists a real number $c_{\widetilde{W}} > \mathfrak{H}(Y)$ such that $f(\Theta_{\widetilde{W}}(\Gamma, \Delta)) \geq c_{\widetilde{W}}$ for all $\Gamma, \Delta \in \Sigma^+(A)$.  Any such real number $c_{\widetilde{W}}$ is referred to as a {\em threshold.}

\begin{lemm}\label{lemm:UnifSeparate} Let $\mS_{\widetilde{W}}$ be uniformly separated above $\Sigma(Y)$ and $c_{\widetilde{W}}$ be a threshold.  Then there exists an integer $R\geq1$ such that, for all $\Theta \in \mS_{\widetilde{W}}$ and all $r\in \ZZ$ with $|r|>R$, we have \[f(\sigma^r \Theta) < c_{\widetilde{W}}.\]  Hence, \[\mathfrak{m}(\Theta) = \max_{|r| \leq R}f(\sigma^r \Theta).\]
 
\end{lemm}

\begin{proof}  Let $c:= c_{\widetilde{W}}$.  Assume that no such $R$ exists.  Then, for every  $k \geq 1$, there exists $\Gamma_k, \Delta_k \in \Sigma^+(A)$ and $r_k \in \ZZ$ with $|r_k| >k$, such that we have \begin{align}\label{eqn:lemm:UnifSeparate}
 f(\sigma^{r_k}\Theta(\Gamma_k, \Delta_k)) \geq c.  \end{align} Hence, we have that \[d_k := \left\lfloor \frac {\min(|r_k + \widetilde{n}_-|,|r_k - \widetilde{n}_+|)} 2\right\rfloor \rightarrow \infty\] as $k \rightarrow \infty$, from which it follows that $[r_k - d_k,   r_k + d_k] \cap [-\widetilde{n}_-, \widetilde{n}_+] = \emptyset$ for all $k$ large enough.  Let $k$ be large enough.  Consequently, $\left(\sigma^{r_k} \Theta(\Gamma_k, \Delta_k)\right) |_{[-d_k, d_k]}$ lies entirely in a concatenation of elements of $A$.  Since the elements of $A$ are finite concatenations of the elements of $Y$, $\left(\sigma^{r_k} \Theta(\Gamma_k, \Delta_k)\right) |_{[-d_k, d_k]}$ is a word of some $\gamma_k \in \Sigma(Y)$ such that \[\gamma_k |_{[-d_k, d_k]} = \left(\sigma^{r_k} \Theta(\Gamma_k, \Delta_k)\right) |_{[-d_k, d_k]}.\]  Since $\Sigma(Y)$ is compact, there is subsequence on which $\gamma_k \rightarrow \theta \in \Sigma(Y)$.  Since $f$ is continuous, we have that $f(\theta) \geq c$ by (\ref{eqn:lemm:UnifSeparate}).  Since $\theta \in\Sigma(Y)$ and $c$ is a threshold, we have that $f(\theta) \leq \mathfrak{H}(Y)< c$, which is a contradiction and gives the desired result.
 
\end{proof}

Let $R$ be as in Lemma~\ref{lemm:UnifSeparate}.  Choose $m_-, m_+ \in \NN \cup \{0\}$ and elements $\delta'_1,\cdots, \delta'_{m_-}, \gamma'_1, \cdots, \gamma'_{m_+} \in A$ such that $n_-:= \widetilde{n}_-+ \sum_{i=1}^{m_-} |\delta'_i|_{\NN} \geq R$ and $n_+:= \widetilde{n}_++ \sum_{i=1}^{m_+} |\gamma'_i|_{\NN} \geq R$ and define the fixed word \[W:= \delta'_{m_-} \cdots \delta'_1\widetilde{W}\gamma'_1 \cdots \gamma'_{m_+},\] which we call the {\em distinguished word}.  After extending the indexing of $\widetilde{W}$ to $W$, we have that $W = a_{-n_-} \cdots a_{-1} a_0 a_1 \cdots  a_{n_+}$.  We call $a_0$ the {\em distinguished digit (of $W$)}.  If $\mS_{\widetilde{W}}$ is uniformly separated above $\Sigma(Y)$, then the subset $\mS_W:= \{\Theta_W(\Gamma, \Delta): \Gamma, \Delta \in \Sigma^+(A)\} \subset \mS_{\widetilde{W}}$ is also uniformly separated above $\Sigma(Y)$.  Every threshold for $\mS_{\widetilde{W}}$ is also a threshold for $\mS_W$.  Note that $a_0$ is at  the $0$-th position of every element of $\mS_W$.  

Let $K^+:= K(A)$ and $K^-:=K(A^\mathsf{T})$.  Define the mapping $\xi: \Sigma^+(A) \rightarrow K^+$ by $\xi(\Gamma) := [0; \gamma_1,\gamma_2, \cdots]$ and the mapping $\eta: \Sigma^+(A) \rightarrow K^-$ by $\eta(\Delta) = [0;\delta_1^\mathsf{T}, \delta_2^\mathsf{T}, \cdots].$  Since $A$ and $A^\mathsf{T}$ are prefix-free, $\xi$ and $\eta$ are, respectively, bijections.

%\begin{lemm}\label{lemm:SWMarkov} For every integer $|r| \leq R$, there exist an $a_r \in \NN$ and $h^+_r, h^-_r  \in \PGL(2, \ZZ)$, defined on open sets containing $K^+$ and $K^-$, respectively, such that, for all $\Gamma, \Delta \in \Sigma^+(A)$, we have \[f(\sigma^r \Theta_W(\Gamma, \Delta))= a_r + h^+_r(\xi(\Gamma))+ h^-_r(\eta(\Delta)).\]
% 
%\end{lemm}

\begin{lemm}\label{lemm:SWMarkov} For every integer $|r| \leq R$, there exist an $a_r \in \NN$ and $h^+_r, h^-_r  \in \PGL(2, \ZZ)$, defined on open sets containing $K^+$ and $K^-$, respectively, such that, for all $\Theta_W(\Gamma, \Delta) \in \mS_W$, we have \[f(\sigma^r \Theta_W(\Gamma, \Delta))= a_r + h^+_r(\xi(\Gamma))+ h^-_r(\eta(\Delta)).\]
 
\end{lemm}

\begin{proof}  Let $\Gamma := \gamma_1 \gamma_2 \cdots$ and $\Delta:= \delta_1\delta_2 \cdots$.  Since $[-R, R ] \subset [-n_-, n_+]$, $a_r$ is a digit of $W$, and thus we can write $W = W_r^- a_r W_r^+$.  Consequently, we have that $f(\sigma^r \Theta_W(\Gamma, \Delta)= a_r + [0;W_r^+,\gamma_1,\gamma_2, \cdots] + [0;(W_r^-)^\mathsf{T}, \delta_1^\mathsf{T},\delta_2^\mathsf{T}, \cdots] )$.  Now if $W_r^+$ is the empty word, then set $h_r^+:= I$ where $I$ is the identity element of $\PGL(2,\ZZ)$.  Similarly for $h_r^-$.  Otherwise, recalling that we identify an element of $\GL(2, \ZZ)$ with its image under the canonical projection map, set \[h_r^+:= M_{W_r^+}, \quad h_r^- := M_{(W_r^-)^\mathsf{T}}.\]   Using elementary properties of continued fractions, it follows that $h_r^+, h_r^-$ have no poles in $[0,1]$.  Since $K^+, K^- \subset [0,1]$, the result follows. \end{proof}

For each $r \in \ZZ$ with $|r|\leq R$, let $a_r$, $h^+_r$, and $h^-_r$ be as in Lemma~\ref{lemm:SWMarkov}.  Define $g_r: K^+ \times K^- \rightarrow \RR$ by $g_r(\xi, \eta):= a_r + h^+_r(\xi) + h^-_r(\eta)$, and let \[G(\xi, \eta) := \max_{|r|\leq R} g_r(\xi, \eta).\]  

\begin{lemm}\label{lemm:BaireChoice}  There exists an integer $|r_*|\leq R$ and $\varphi_*^+, \varphi_*^-  \in \PGL(2, \ZZ)$, defined on open sets containing $K^+$ and $K^-$, such that \begin{align*}
&C^+:= \varphi_*^+(K^+) \subset K^+, \quad C^-:= \varphi_*^-(K^-)\subset K^-,\\ &\dim_H C^+ = \dim_H K^+, \quad \dim_H C^- = \dim_H K^-\end{align*} and, for all $(\xi, \eta) \in C^+ \times C^-$, we have $G(\xi, \eta) = g_{r_*}(\xi, \eta)$.
\end{lemm}

\begin{proof} For $|A|=1$, $K^+$ and $K^-$ are singleton sets.  Choosing $r_*$ such that $g_{r_*}(\xi, \eta) \geq g_r(\xi, \eta)$ for all integers $|r| \leq R$ and the unique element $(\xi, \eta)$ of $K^+ \times K^-$ and setting $\varphi_*^+ = I=\varphi_*^- $, where $I$ is the identity in $\PGL(2, \ZZ)$, yields the desired result.

Hence, we may assume $|A| \geq 2$.  Then $K^+$ and $K^-$ are regular Cantor sets and thus compact.  For each integer $r$ with $|r| \leq R$, define \[D_r:= \{(\xi, \eta) \in K^+ \times K^-:  g_r(\xi, \eta) = G(\xi, \eta)\}.\]  Since $g_r$ and $G$ are continuous, the sets $D_r$ are closed subsets of the compact metric space $K^+ \times K^-$.  By the Baire Category Theorem, there exists some $D_{r_*}$ such that $\inte (D_{r_*})$ is a nonempty open subset of $K^+ \times K^-$.  Consequently, there are nonempty open sets $U^+ \subset K^+$ and $U^- \subset K^-$ such that $U^+ \times U^- \subset D_{r_*}$.  If $\xi \in U^+$, then $\xi = [0; \gamma_1, \gamma_2, \cdots]$ for $\gamma_n \in A$.  For each $n \in \NN$, let $\Gamma_n:= \gamma_1 \cdots \gamma_n$.  Now, for all $n \in \NN$, we have that \[\xi \in M_{\Gamma_n}(K^+) \subset K^+.\]  Since $\diam (M_{\Gamma_n}(K^+)) \rightarrow 0$ as $n \rightarrow \infty$, we have that there exists a large enough integer $n$ such that $M_{\Gamma_{n}}(K^+) \subset U^+$.   Let $\varphi^+_* :=M_{\Gamma_{n}}$ for this $n$.  Now there are no poles of $\varphi^+_*$ in $[0,1]$, and thus $\varphi^+_*$ is a $C^\infty$-diffeomorphism on an open set containing $K^+$, which gives $\dim_H(C^+) = \dim_H(K^+)$.   Similarly, for $U^-$, we have $\varphi^-_*$ and the dimension result.  This yields the desired result. \end{proof}

Let $r_*$, $C^+$, $C^-$,$\varphi_*^+$, and $\varphi_*^-$ be as in Lemma~\ref{lemm:BaireChoice}.  Define \begin{align}\label{eqn:defnSumSet}
 V:= G(C^+\times C^-) = g_{r_*}(C^+\times C^-) = a_{r_*}+h^+_{r_*}(C^+)+ h^-_{r_*}(C^-). \end{align}

%Applying Lemma~\ref{lemm:UnifSeparate}, we note that, for all $\lambda \in V$, there exists $(\xi, \eta) \in C^+ \times C^-$, corresponding to $\xi :=\xi(\Gamma)$ and $\eta:=\eta(\Delta)$, such that $\lambda = \mathfrak{m}(\Theta_W(\Gamma, \Delta))\geq c_W$ where $c_W$ is a threshold for $\mS_W$.

\begin{lemm}\label{lemm:CantorSetSumMarkovValues}
For all $\lambda \in V$, there exists $\Gamma, \Delta \in \Sigma^+(A)$ such that $\lambda =  f(\sigma^{r_*} \Theta_W(\Gamma, \Delta))= \mathfrak{m}(\Theta_W(\Gamma, \Delta))\geq c_W$ where $c_W$ is a threshold for $\mS_W$.  Moreover, $V$ is compact and \[\dim_H V = \min(1, 2 \dim_H K(A)).\]
\end{lemm} 

\begin{proof}  Let $\lambda = G(\xi, \eta)$.  There exists $ \Gamma, \Delta \in \Sigma^+(A)$ such that $\xi= \xi(\Gamma)$ and $\eta = \eta(\Delta)$. Then Lemmas~\ref{lemm:UnifSeparate} and~\ref{lemm:SWMarkov} give $\lambda = G(\xi, \eta)= f(\sigma^{r_*} \Theta_W(\Gamma, \Delta))=\mathfrak{m}(\Theta_W(\Gamma, \Delta)) \geq c_W$ where $c_W$ is a threshold for $\mS_W$.

Compactness of $V$ follows from Lemma~\ref{lemm:BaireChoice}, the continuity of $g_{r_*}$, and the compactness of $K^+\times K^-$.  If $|A|=1$, then $K^+$ and $K^-$ are singleton sets, and, by Lemma~\ref{lemm:BaireChoice}, $V$ is also a singleton set.  This yields the desired dimension formula.  

Otherwise, we have that $|A|\geq 2$.  Since $h^+_{r_*}\circ\varphi_*^+$ is a $C^\infty$-diffeomorphism on an open set containing $K^+$ and similarly for $K^-$, Lemma~\ref{lemm:Moreira} applies to $h^+_{r_*}(C^+)$ and $h^-_{r_*}(C^-)$.  (In fact, both are non-essentially affine.)  Now Lemma~\ref{lemm:BaireChoice} and (\ref{eqnHDKHDKTrans}) yield the desired result.
\end{proof}

\subsection{Projective cocycle steering}\label{subsec:ProjCocycleSteer}  In this section, we return to our setup over $R_Q$.  We retain the notation from Section~\ref{sec:DistWord}.  Let $H:=H_Q(Y)$ and $J \subset R_Q$ be nonempty.  In addition, we require that $\rho_Q(W) \in H$ and \begin{align}\label{eqn:ProjsteeringAisIdentity}
 \rho_Q(a) = 1_H \end{align} for all $a \in A$.  Write $W= W_{r_*}^- a_{r_*}W_{r_*}^+$, and define $M_* := \rho_Q(W_{r_*}^-) \in \PGL(2, R_Q)$.   Note that $W_{r_*}^-$ may be the empty word, and, similarly, $W_{r_*}^+$ may also be.  Recall, we are using the convention that $\rho_Q$ applied to the empty word is $1_H$.  We say $(W, r_*, H)$ is {\em $J$-accessible} if, for every $j \in J$, there exists $h_j \in H$ such that \[h_jM_*[0:1] = [-j:j+1] \quad \textrm{ in } \PP^1(R_Q).\]  Fix $\lambda \in V$. Let $\Gamma:=\gamma_1\gamma_2\cdots, \Delta:=\delta_1 \delta_2 \cdots \in \Sigma^+(A)$ be as in Lemma~\ref{lemm:CantorSetSumMarkovValues}.  For integers $m\geq 1$, set \[\tau_m:= \delta_m \delta_{m-1} \cdots\delta_1 W \gamma_1 \gamma_2 \cdots \gamma_m,\] where the $W$ is referred to as the {\em designated $W$ (of $\tau_m$)}.  Note that $\rho_Q(\tau_m) = \rho_Q(W)$ for every $m$. Let  $\mathfrak{L}^+_m:=|\gamma_1 \gamma_2 \cdots \gamma_m|_\NN$ and $\mathfrak{L}^-_m:=|\delta_m \delta_{m-1} \cdots\delta_1|_\NN$.

Assume $(W, r_*, H)$ is $J$-accessible (we prove this for the applications in the form of transitivity).  Then there exists an $s_j, \widetilde{s}_j \in Y^+$ such that $\rho_Q(s_j) = h_j$ and $\rho_Q(\widetilde{s}_j) = (h_j \rho_Q(W))^{-1}$ by construction of $H$.  Define the word \[\Pi_{m,j} := s_j \tau_m \widetilde{s}_j\] where the designated $W$ of $\tau_m$ is also referred to as the {\em designated $W$ (of $\Pi_{m,j}$)}.  Note that \begin{align}\label{eqn:PimjGivesIdentity}
 \rho_Q(\Pi_{m,j})= 1_H \end{align} for every $m$ and $j$.  Since we consider the digits of $\Pi_{m,j}$ from left to right, we note that, roughly speaking, $s_j$ steers the cocycle $\rho_Q(\delta_m \delta_{m-1} \cdots\delta_1W_{r_*}^-)$ so that it takes $[0:1]$ to $[-j:j+1]$ and $\widetilde{s}_j$ resets the cumulative cocycle $\rho_Q(s_j \tau_m)$ back to $1_H$.  

% \[n_{m,j}:= |\Xi_{m,j}|_\NN + |s_j|_\NN + \sum_{i=1}^m|\delta_i|_\NN + |W_{r_*}^-|_\NN.\]  

Suppose we have that $\Theta_{m,j} := \Xi_{m,j} \Pi_{m,j} \widetilde{\Xi}_{m,j}$ where $\Xi_{m,j}$ is the empty word or a word over $\NN$ and $\widetilde{\Xi}_{m,j}$ is the empty word, a  word over $\NN$, or an element of $\NN^\NN$.  If $\Theta_{m,j} \in \NN^\NN$, write $\Theta_{m,j} = \theta_1 \theta_2 \cdots.$ Otherwise there exists $k \in \NN$ such that $\Theta_{m,j} = \theta_1 \theta_2 \cdots \theta_k \in \NN^k.$  Define \begin{align}\label{eqn:posnmjinThetamj}
 n_{m,j}:= |\Xi_{m,j}|_\NN + |s_j|_\NN +\mathfrak{L}^-_m + |W_{r_*}^-|_\NN. \end{align}  Consequently, \begin{align}\label{eqn:convergentposition}
 \theta_{n_{m,j}+1} = a_{r_*}. \end{align}

% $(W, r_*, H)$ is $J$-accessible 
\begin{lemm}\label{lemm:accessibleConvergCond}  Assume $\rho_Q(\Xi_{m,j}) =1_H$, and let \[\frac{p_{n_{m,j}}}{q_{n_{m,j}}} = [0; \theta_1, \cdots, \theta_{n_{m,j}}].\]  Then $[p_{n_{m,j}}:q_{n_{m,j}}] = [-j:j+1]$ in $\PP^1(R_Q)$.
 
\end{lemm}

\begin{proof}
 Fix $m\geq 1$ and $j \in J$.  Let $\ell:= n_{m,j}$.  Thus, by (\ref{eqn:ProjsteeringAisIdentity}), we have  $\rho_Q(\theta_1 \cdots \theta_\ell) = h_j M_*$ .  By (\ref{eqn:ContFracMatConvergents}), we $[p_\ell:q_\ell] = \rho_Q(\theta_1 \cdots \theta_\ell)[0:1] =  h_j M_* [0:1] = [-j: j+1]$.\end{proof}

Now enumerate $J = \{j_1, \cdots j_s\}$ and let \[\Pi_\lambda:= \Pi_{1,j_1}\cdots\Pi_{1,j_s}\Pi_{2,j_1} \cdots \Pi_{2,j_s} \Pi_{3,j_1} \cdots.\]  Note that $\Pi_\lambda \in \NN^\NN$ and thus we may write $\Pi_\lambda = \theta_1 \theta_2 \cdots$ for $\theta_i \in \NN$.  Let $\mu_{m,j}\in \NN$ denote the position of the digit in $\Pi_\lambda$ of the distinguished digit of the designated $W$ of $\Pi_{m,j}$ and let $\mathcal{M}(\Pi_\lambda) := \{\mu_{m,j} : m \geq 1, j \in J\} \subset \NN$.  Note that the set $\mathcal{M}(\Pi_\lambda)$ becomes sparser as the elements get larger.\footnote{Our proofs remain unchanged should $\Pi_\lambda$ have subwords $W$ that are not designated.}  Now define the continued fraction \[y_\lambda := [0; \theta_1, \theta_2, \cdots ]\] and let $p_n/q_n$ be its convergents.  Let $n_k \rightarrow \infty$ be a sequence of natural numbers.  We say {\em $\Pi_\lambda$ converges locally to $b:=(b)_{\ell\in \ZZ} \in \Sigma$ along $(n_k)_{k \in \NN}$} if for every integer $L \geq 1$, there exists $k_L \in \NN$ such that, for all $k \geq k_L$ and all integers $|\ell|\leq L$ ,we have that $\theta_{n_k+1 +\ell} = b_\ell$.

\begin{lemm}\label{lemm:PiLambdaLocConverg}  Let $n_k \rightarrow \infty$ be a sequence of natural numbers.  Suppose $\Pi_\lambda$ converges locally to $b \in \Sigma$ along $(n_k)_{k \in \NN}$.  Then either $b = \sigma^t \Theta_W(\Gamma, \Delta)$ for some $t \in \ZZ$ or $b \in \Sigma(Y)$.
 
\end{lemm}
\begin{proof}  Let $d_k:= \min\{|n_k + 1  - \mu| : \mu \in \mathcal{M}(\Pi_\lambda)\}$.   After passing to a subsequence if needed, we may assume, without loss of generality, that either $\sup d_k < \infty$ or  $d_k \rightarrow \infty$ hold.  Consider Case 1:   $\sup d_k < \infty$.  For every $k \in \NN$, choose $\mu_k:=\mu_{m_k,j_k} \in  \mathcal{M}(\Pi_\lambda)$ such that $d_k = |n_k + 1  - \mu_k|$ and let $t_k = n_k + 1  - \mu_k$.  Since $t_k \in [-\sup d_k, \sup d_k ]$, we may, after passing to a subsequence if needed, assume, without loss of generality, that there exists a $t \in \ZZ$ such that $t_k = t$ for all $k \in \NN$.  Since $n_k \rightarrow \infty$ but $t$ is fixed, we have that $m_k \rightarrow \infty$.  Thus, we have $\mathfrak{L}^+_{m_k} \rightarrow \infty$ and $\mathfrak{L}^-_{m_k} \rightarrow \infty$.  Fix an integer $L \geq 1$.  Consequently, for all sufficiently large $k$, we have \[\Pi_\lambda |_{[\mu_k+t -L, \mu_k+t +L]}\] is a subword of $\tau_{m_k}$ that is itself a subword of $\Pi_{m_k, j_k}$.  Recall that $\mu_k$ is the position in $\Pi_\lambda$ of the distinguished digit in the subword $\Pi_{m_k, j_k}$. Hence, we have that \[b_\ell=\theta_{n_k+1+\ell}=\theta_{\mu_k+t +\ell}= \left(\sigma^t \Theta_W(\Gamma, \Delta) \right)_\ell\] for every $|\ell| \leq L$ and all sufficiently large $k$. Here, the first equality follows by local convergence.  Since $L$ is arbitrary, we have that $b = \sigma^t \Theta_W(\Gamma, \Delta)$.

Now consider Case 2:  $d_k \rightarrow \infty$.  Fix an integer $L \geq 1$.  Recall that we have $W = a_{-n_-} \cdots a_{-1} a_0 a_1 \cdots  a_{n_+}$, and set $I_W := [-n_-, n_+] \cap \ZZ$ and $r_W:= \max(n_-, n_+)$.  Consequently, we have that $d_k > L + r_W$ for all $k$ sufficiently large.  Therefore, for every $p \in [n_k+1 -L, n_k+1 +L] \cap \ZZ$ and every $\mu \in \mathcal{M}(\Pi_\lambda)$, we have that \begin{align}\label{eq1:lemm:PiLambdaLocConverg}
 |p - \mu| \geq |n_k+1 -\mu| - |n_k+1 -p| \geq d_k - L >r_W. \end{align}  This implies that $p \notin \mu+I_W$ and, consequently, \begin{align}\label{eq2:lemm:PiLambdaLocConverg}
 [n_k+1-L, n_k+1+L] \cap \bigcup_{\mu \in \mathcal{M}(\Pi_\lambda)} \mu + I_W = \emptyset.  \end{align} 
  
Now setting $L_k := \lfloor (d_k - r_W)/2\rfloor$ for $L$ in (\ref{eq1:lemm:PiLambdaLocConverg}, \ref{eq2:lemm:PiLambdaLocConverg}) implies that \[\theta |_{[n_k+1-L_k, n_k+1+L_k]}\] cannot overlap with any designated $W$ for all sufficiently large $k$.  Since, aside from the designated $W$'s, the words concatenated to form $\Pi_\lambda$ come from $Y$, we may choose $w^{(k)} \in \Sigma(Y)$ such that $w^{(k)}_\ell = \theta_{n_k+1+\ell}$ for integers $|\ell|\leq L_k$.  Fix an integer $R \geq1$.  Then for all $k$ sufficiently large, we have that $R \leq L_k$ and, thus, local convergence gives that $w^{(k)}_\ell = b_\ell$ for all integers $|\ell|\leq R$.  Since $\Sigma(Y)$ is compact, we have that $b \in \Sigma(Y)$. \end{proof}

\begin{lemm}\label{lemm:GoodconvergentsLagrange}  Let $p_n/q_n$ denote the $n$th-order convergent of the irrational number $y_\lambda$.  Then we have that \[k(y_\lambda) = \lambda.\]  Moreover, for every $j \in J$, there exists a sequence of natural numbers $n_{m,j}\rightarrow \infty$ such that \[[p_{n_{m,j}}:q_{n_{m,j}}] = [-j:j+1] \textrm{ in } \PP^1(R_Q) \textrm{ for all } m \geq 1 \quad \textrm{ and} \quad \lim_{m \rightarrow \infty} \Lambda_{n_{m,j}}(y_\lambda) = \lambda.\]
 
\end{lemm}
\begin{proof}
By Lemma~\ref{lemm:CantorSetSumMarkovValues}, we have that there exists $\Gamma, \Delta \in \Sigma^+(A)$ such that \[\lambda = f(\sigma^{r_*} \Theta_W(\Gamma, \Delta))=\mathfrak{m}(\Theta_W(\Gamma, \Delta))\geq c_W > \mathfrak{H}(Y)\] where the final inequality follows because $\mS_W$ is uniformly separated above $\Sigma(Y)$.  Now let  $\Xi_{m,j}$ denote the prefix of $\Pi_\lambda$ before the $\Pi_{m,j}$ and $\widetilde{\Xi}_{m,j}$ denote the suffix after the $\Pi_{m,j}$, giving $\Pi_\lambda = \Xi_{m,j}\Pi_{m,j}\widetilde{\Xi}_{m,j}$.  Consequently, we have the $\rho_Q(\Xi_{m,j})=1_H$ by (\ref{eqn:PimjGivesIdentity}).  Defining $n_{m,j}$ as in (\ref{eqn:posnmjinThetamj}) yields (\ref{eqn:convergentposition}) and $n_{m,j}+1 = \mu_{m,j}+ r_*$.  We note that $\Pi_\lambda$ converges locally to $\sigma^{r_*} \Theta_W(\Gamma, \Delta)$ along the sequence $(n_{m,j})$, and the proof is identical to that of Case 1 of Lemma~\ref{lemm:PiLambdaLocConverg} with $t$ replaced by $r_*$.  Consequently, (\ref{eqnPerronSubseqConv}) gives the desired limit and (\ref{eqnPerronConv}) yields $k(y_\lambda) \geq \lambda$.  The equality in $\PP^1(R_Q)$ now follows by Lemma~\ref{lemm:accessibleConvergCond}.

By construction, the partial quotients of $y_\lambda$ are bounded and thus $k(y_\lambda)< \infty$.  Let $D$ be the finite set of partial quotients of $y_\lambda$. Now choose a sequence of natural numbers $n_k \rightarrow \infty$ such that $\Lambda_{n_k}(y_\lambda) \rightarrow k(y_\lambda)$.  Extend $\Pi_\lambda$ to the left with elements of $D$ to obtain an element $\widetilde{\Pi}_\lambda$ of $\Sigma(D)$.  Consider the sequence $\sigma^{n_k+1}\widetilde{\Pi}_\lambda$ in $\Sigma(D)$.  Since $\Sigma(D)$ is compact, $\sigma^{n_k+1}\widetilde{\Pi}_\lambda$ has a convergent subsequence $\sigma^{n_{k_i}+1}\widetilde{\Pi}_\lambda$ to some $b \in \Sigma(D)$.  Consequently, we have $\Pi_\lambda$ converges locally to $b$ along $(n_{k_i})$.  By (\ref{eqnPerronSubseqConv}), we have that $k(y_\lambda) = f(b)$.  Now, by Lemma~\ref{lemm:PiLambdaLocConverg}, either $b = \sigma^t \Theta_W(\Gamma, \Delta)$ for some $t \in \ZZ$ or $b \in \Sigma(Y)$.  Thus either $k(y_\lambda) = f(\sigma^t \Theta_W(\Gamma, \Delta)) \leq \mathfrak{m}(\Theta_W(\Gamma, \Delta)) = \lambda$ or $k(y_\lambda) \leq \mathfrak{H}(Y) < c_W \leq \lambda$.  Thus, we have $k(y_\lambda) = \lambda$.

\end{proof}

\begin{prop}\label{prop:projcocyclesteering}  Let $0 \in J \subset R_Q$, $H:= H_Q(Y)$, $\mS_W$ be uniformly separated above $\Sigma(Y)$ with threshold $c_W$, $(W, r_*, H)$ be $J$-accessible, $\rho_Q(W) \in H$, and $\rho_Q(a) = 1_H$ for every $a \in A$.  Then $QV$ is a compact set such that \[QV \subset \eL\left(\bigcap_{j \in J}E_{j/Q} \right) \quad \textrm { and } \quad \dim_H(QV) = \min(1, 2 \dim_H K(A)).\]  Moreover, if $s \in \RR$ satisfies $\mathfrak{m}(\Theta_W(\Gamma, \Delta)) < s$ for every $\Gamma, \Delta \in \Sigma^+(A)$, then $QV \subset (-\infty, Qs)$.
 
\end{prop}

\begin{proof} Fix $\lambda \in V$.  By Lemma~\ref{lemm:GoodconvergentsLagrange} and Hurwitz's theorem, we have that $k(y_\lambda) = \lambda \geq \sqrt 5$ and that the condition of Lemma~\ref{lemm:LagrangeUnderQPush} is satisfied, thus, giving $k(A_j y_\lambda) = Q \lambda$ for all $j \in J$.  Setting $x_\lambda := A_0(y_\lambda)$ yields that \[k(x_\lambda) =Q\lambda= k(x_\lambda+ j/Q)\] for all $j \in J$.  Since multiplication by $Q$ is bilipschitz, Lemma~\ref{lemm:CantorSetSumMarkovValues} yields the dimension result.  Lemma~\ref{lemm:CantorSetSumMarkovValues} also yields the result involving $s$.
\end{proof}

%Therefore, the local convergence misses every copy of $W$ in every $\Pi_{m, j}$ and thus $b \in \Sigma(Y)$.

\section{Transitive projective cocycle steering and applications}\label{sec:TransitiveSteering} In this section, we establish transitive projective cocycle steering (Theorem~\ref{prop:transitivesteering}), giving lower bounds when $H_Q$ acts transitively on $\PP^1(R_Q)$, from which follows a full Hausdorff dimension regime (Corollary~\ref{cor:FullHD}) and a positive Hausdorff dimension regime (Proposition~\ref{prop:PosHD})  for $d_F$ and a scaling law (Proposition~\ref{prop:ScalingLaw}).

\subsection{Proof of Theorem~\ref{prop:transitivesteering}}\label{subsubsec:prooftransitivesteering} Let $d := \dim_HK(X)$ and $H := H_Q(X)$.  If $d = 0$, the conclusion is immediate.  Otherwise, we have $d> 0$.  Fix $0 < \varepsilon < d/2$.   Applying Lemma~\ref{lemm:HighXRemoveExtremes} with $R=2$ gives a natural number $n$ such that, for $B:= X^n$, $B_{\ext}:= \{\min(B), \max(B)\}$, and $Y:= B \backslash B_{\ext}$, we have $H_Q(Y) = H$ and $\dim_H K(Y) > d - \varepsilon$.  Note that the minimum and maximum is with respect to our strict total ordering on $B$.  Now note that we have \[\mathfrak{H}(B) = \sup_{\theta \in \Sigma(B)} \sup_{n \in \ZZ} f(\sigma^n \theta) = \sup_{\theta \in \Sigma(B)} \mathfrak{m}(\theta) = \sup_{\theta \in \Sigma(X)} \mathfrak{m}(\theta)= \mathfrak{H}(X) < s\] by $\Sigma(B) = \Sigma(X)$ and $\sigma$-invariance.  By the argument in the proof of~\cite[Lemma~3]{Mor18}, we also have that $\mathfrak{H}(Y) < \mathfrak{H}(B)$.  Applying Lemma~\ref{lemm:IdentityAlphabet} to $Y$ yields a finite prefix-free alphabet $A$ consisting of concatenations of the elements of $Y$ such that $A^\mathsf{T}$ is also prefix-free.  Thus, Sections~\ref{sec:DistWord} and~\ref{subsec:ProjCocycleSteer} apply.  Also, we have that $\rho_Q(a) = 1_H$ for every $a \in A$.  Let $J=\{0\} \cup \{Qr \mod Q: r \in F\}$.

By continuity of $f$ and compactness of $\Sigma(B)$, choose $\beta \in \Sigma(B)$ such that $f(\beta) = \mathfrak{H}(B)$.  Let $\beta = (b_j)_{j \in \ZZ}$ where $b_j \in \NN$, and, thus, we have $\mathfrak{H}(B) = b^+ + b^-$ where $b^+:=[b_0; b_1, \cdots]$ and $b^-:= [0; b_{-1}, b_{-2}, \cdots]$.  We also have that $\beta = \cdots \beta_{-2} \beta_{-1} \beta_0 \beta_1 \beta _2 \cdots$ where $\beta_i \in B$ and $b_0$ lies in $\beta_0$. For every integer $N \geq 1$, define \begin{align}\label{eqn:WTildeNDefn}
 \widetilde{W}_N := \beta_{-N} \beta_{-N+1} \cdots \beta_{-1} \beta_0 \beta_1 \beta_2 \cdots \beta_N \end{align} and let $b_0$ be its distinguished digit.

\begin{lemm}\label{lemm:transSWthres}
For all sufficiently large $N \in \NN$, we have that $\mS_{\widetilde{W}_N}$ is uniformly separated above $\Sigma(Y)$ with threshold $c:= c_{\widetilde{W}_N} = (3 \mathfrak{H}(B) + \mathfrak{H}(Y))/4 > \mathfrak{H}(Y)$.
\end{lemm}

\begin{proof}
 Let $\theta \in \mS_{\widetilde{W}_N}$ with $0$-th position at the distinguished digit of $\widetilde{W}$ and write $\theta:=(\theta_j)_{j \in \ZZ}$ where $\theta_j \in \NN$.  Thus $\theta_0 = b_0$.  Let $\theta^+:= [\theta_0; \theta_1, \theta_2, \cdots]$ and $\theta^-:= [0; \theta_{-1}, \theta_{-2}, \cdots]$.  Hence, we have that $f(\theta) = \theta^+ + \theta^-$.  Now $b^+$ and $\theta^+$ agree on at least the first $N$ digits and, likewise, $b^-$ and $\theta^-$ agree on at least the first $N$ digits. Using elementary properties of continued fractions, we have $|b^+ - \theta^+| \leq F_{N+1}^{-2}$ and $|b^- - \theta^-| \leq F_{N+1}^{-2}$, where $F_N$ is the $N$-th term of the Fibonacci sequence (with $F_1 = 1 = F_2$).  Thus we have that \[|f(\theta) - \mathfrak{H}(B)| \leq 2F_{N+1}^{-2} \leq (\mathfrak{H}(B) - \mathfrak{H}(Y))/4\] for all sufficiently large $N$, yielding the desired result. \end{proof}

Choose $N$ large enough to satisfy Lemma~\ref{lemm:transSWthres}.  Using $\widetilde{W} = \widetilde{W}_N$, form the distinguished word $W$.  Therefore, $\mS_W \subset \Sigma(B) = \Sigma(X)$.  Consequently, for every $\Theta \in \mS_W$, we have $\mathfrak{m}(\Theta) <s$.  Now recall that $Y$ depends on $\varepsilon$.  Let $V_\varepsilon$ be the set $V$ in (\ref{eqn:defnSumSet}) for this $Y$.  Then Lemmas~\ref{lemm:CantorSetSumMarkovValues}and~\ref{lemm:IdentityAlphabet}  give $\dim_H V_\varepsilon = \min(1, 2 \dim_H K(A)) \geq \min(1, 2(d - 2 \varepsilon))$.

Finally, since $H$ acts transitively on $\PP^1(R_Q)$, it follows immediately that $(W, r_*,H)$ is $J$-accessible.  Applying Proposition~\ref{prop:projcocyclesteering} yields \[QV_\varepsilon \subset \eL\left(\bigcap_{j \in J} E_{j/Q} \right) \cap (-\infty, Qs),\] and, recalling that $k(x+\ZZ) = k(x)$, we have $d_F(Qs) \geq \min(1, 2(d - 2 \varepsilon))$.  Letting $\varepsilon \rightarrow 0$ yields the desired result and concludes the proof of Theorem~\ref{prop:transitivesteering}.

\subsection{A compact set of full Hausdorff dimension}\label{subsec:FullHD}  In this section, we prove one of the main results of this paper Theorem~\ref{thm:MainBad}~\ref{it7:FullDimCpctSet:MainBad} using transitive projective cocycle steering.  Arithmetic enters strongly again in the form of strong approximation, which gives (\ref{eqn:StrongApprForAllModulus}).  First, we choose a suitable alphabet:  \[X:= \{11, 12, 21, 22\}.\]  This is the alphabet considered by Bourgain and Kontorovich~\cite[(1.32) and Remark~1.33]{BK14}.  After reduction modulo $Q$, strong approximation~\cite[Theorem in Section~0]{MVW84} gives, for all sufficiently large primes $Q$, \begin{align}\label{eqn:StrongApprForAllModulus} \langle M_w: w \in X\rangle = \SL(2, R_Q), 
  \end{align} which acts transitively on $\PP^1(R_Q)$.  For $X$, (\ref{eqn:StrongApprForAllModulus}) actually holds for all integers $Q \geq 2$~\cite[Remark~1.33]{BK14}.  Taking images under the canonical projection map, we have that $H:=H_Q(X) = \PSL(2, R_Q) \trianglelefteq \PGL(2, R_Q)$ acts transitively on $\PP^1(R_Q)$.  Note that $\Sigma(X) = \Sigma(\{1, 2\})$ and $\mathfrak{h}(X) = 2 \sqrt 3$~\cite[Page~35]{CF89}. 
  
\begin{lemm}\label{lemm:MaxBinaryHeight}  Let $\theta \in \Sigma(X)$.  Then $f(\theta) = \mathfrak{h}(X)$ if and only if $\theta = \cdots 1212 \cdots$ with $2$ in the $0$th-position.
 
\end{lemm}
\begin{proof}  Let $\theta := (\theta_j)_{j \in \ZZ}$.  If $u := \max K(X)$ and $v:= \min K(X)$, then $M_1(v) = u$ and $M_2(u)=v$.  Thus, we have that $u^2+2u-2=0$ and $u = \sqrt 3 -1 = [0;1, 2, 1, 2 \cdots]$.  Hence, we have that $f(\theta) = \theta_0 + [0;\theta_1, \theta_2, \cdots] + [0;\theta_{-1}, \theta_{-2}, \cdots] \leq 2 + 2 u$.  In order to attain equality, we must have $\theta_0=2$ and $[\theta_1, \theta_2, \cdots] = u = [\theta_{-1}, \theta_{-2}, \cdots].$  Since $<$ is a strict total ordering on $K(X)$, the maximum is unique and thus only $\theta = \cdots 1212 \cdots$ with $2$ in the $0$th-position attains the maximum.

%Hence, with $2$ in the $0$-th position, we have that $f(\cdots 1212 \cdots) = 2 +  2 \cdot [0; 1, 2, 1, 2, \cdots] = 2+ 2 u = 2 \sqrt 3$.
 
\end{proof}  
 \begin{prop}\label{prop:FullHDSet} For every integer $Q\geq2$, there exists a real number $h_Q$ such that  $ 3 < h_Q < 2 \sqrt 3$ and, for every nonempty finite $F \subset \QQ$ with $Qr \in \ZZ$ for all $r \in F$, we have $d_F(t) = 1$ for all $t \geq Qh_Q.$  Moreover, for every $\varepsilon >0$, there exists a compact set $V_{Q, \varepsilon}$ with $\dim_H V_{Q, \varepsilon}=1$ such that \[V_{Q, \varepsilon} \subset \eL(E_F) \cap (2 Q \sqrt 3 - \varepsilon, 2 Q \sqrt 3)\] for every such $F$.
 
\end{prop}
\begin{proof}
Let $\eta := \frac 1 2 (\dim_H K(X) - 1/2)>0$~\cite{Hen89}.  By Lemma~\ref{lemm:HighXRemoveExtremes}, there exists $n \in \NN$ such that, for $Y:= X^n \backslash \{(12)^n, (21)^n\}$, we have $H_Q(Y) = H$ and $\dim_H K(Y) > \dim_H K(X) -\eta > 1/2$.  Now note that every word over $\{1, 2\}$ of length $2n$ with no two adjacent equal digits is either $(12)^n$ or $(21)^n$.  Consequently, every element of $Y$ has either $11$ or $22$ as a subword, which implies that $\mathfrak{h}(Y)< \mathfrak{h}(X)$ by Lemma~\ref{lemm:MaxBinaryHeight}, the continuity of $f$, and the compactness of $\Sigma(Y)$.  Now setting $h_Q := \frac 1 2 \left(\mathfrak{h}(X) + \max(3, \mathfrak{h}(Y))\right)$ and applying Theorem~\ref{prop:transitivesteering} yields the first assertion.

Let $B = X^n$ and $\varepsilon >0$.  Lemma~\ref{lemm:IdentityAlphabet} gives a finite, prefix-free alphabet $A \subset Y^+$ such that $A^\mathsf{T}$ is also prefix-free, $\rho_Q(A) = 1_H$, and $\dim_H K(A) > 1/2$.  Choose $\delta$ such that $0< \delta < \min(\varepsilon/Q, \mathfrak{h}(X)- \mathfrak{h}(Y), \mathfrak{h}(X) - 3)$ and define $c:= \mathfrak{h}(X) -\delta$. Let $\beta = (b_j)_{j \in \ZZ}$ where $b_{2j} = 2$ and $b_{2j+1} = 1$.  We also have $\beta = \cdots \beta_{-2} \beta_{-1} \beta_0 \beta_1 \beta _2 \cdots$ where $\beta_i \in B$ and $b_0$ lies in $\beta_0$. Thus $f(\beta) = \mathfrak{H}(X)$ by Lemma~\ref{lemm:MaxBinaryHeight}.  For every integer $N \geq 1$, define $\widetilde{W}_N$ as in (\ref{eqn:WTildeNDefn}) with $b_0$ its distinguished digit.  By the analog of Lemma~\ref{lemm:transSWthres}, we have that, for large enough $N \in \NN$, $\mS_{\widetilde{W}_N}$ is uniformly separated above $\Sigma(Y)$ with threshold $c >  \mathfrak{h}(Y)$.  

Using $\widetilde{W} = \widetilde{W}_N$, form the distinguished word $W$.  Therefore, $\mS_W \subset \Sigma(B) = \Sigma(X)$. Since $A \subset Y^+$, every $\Theta:=\Theta_W(\Gamma, \Delta) \in \mS_W$ has tails consisting of the elements of $Y$ and thus infinite occurrences of either $11$ or $22$ in both tails.  Thus, by Lemmas~\ref{lemm:UnifSeparate} and ~\ref{lemm:MaxBinaryHeight}, we have $\mathfrak{m}(\Theta) <\mathfrak{h}(X)$ for every $\Theta \in \mS_W$.  Let $V$ be as in (\ref{eqn:defnSumSet}).  Then $V \subset [c, \mathfrak{h}(X)) \subset (\mathfrak{h}(X)- \varepsilon/Q, \mathfrak{h}(X))$.  Let $r^*$ be as in Lemma~\ref{lemm:BaireChoice} and $J:=R_Q$.  Then, by the transitivity of the $H$ action, we have that $(W, r^*, H)$ is $J$-accessible.  Consequently, Proposition~\ref{prop:projcocyclesteering} yields $QV \subset \eL(E_F)$ and $\dim_H(QV)=1$.  Setting $V_{Q, \varepsilon}= QV$ yields the desired result. \end{proof}

%\cite[Theorem in Section~0]{MVW84}, specifically as in~\cite[Remark~1.33]{BK14}.  

\subsection{Positive Hausdorff dimension}  In this section, we prove one of the main results of this paper Theorem~\ref{thm:MainBad}~\ref{it4:PosDim:thm:MainBad} using transitive projective cocycle steering.  Arithmetic enters strongly again in the form of fiber transitivity, Lemma~\ref{lemm:FiberTransivitiy}, and the Chinese Remainder Theorem.  First, we choose a suitable alphabet: \[X_m:= \{u_m, v_m\}, \quad u_m:= 2 1^{2m} 2, \quad v_m= 2 1^{2m+2} 2.\]  Moreira considers this alphabet~\cite[Page~158]{Mor18}, noting that $\mathfrak{H}(X_m) \leq 3+2^{-m}$.  Let \[S := M_{11}, \quad T := M_{22}, \quad U_-:= \begin{pmatrix}  1 & 0 \\ 1 & 1 \end{pmatrix}, \quad U_+:= \begin{pmatrix}  1 & 1 \\ 0 & 1 \end{pmatrix}.\]  Note the importance of $\langle U_+\rangle $ and $\langle U_- \rangle$ is that they act transitively on the affine charts $\{[j:1]: j \in R_Q\}$ and $\{[1:j]: j \in R_Q\}$, respectively.

\begin{lemm}\label{lemm:STConjHQXm}
We have that $H_Q(X_m)$ is conjugate to the subgroup $\langle S, T \rangle$ in $\PGL(2, R_Q)$. 
\end{lemm}
\begin{proof}
 We have that $R:=M_{u_m}^{-1} M_{v_m} = M_2^{-1} S M_2$.  Then $M_{u_m} R^{-m} = T$, and hence $H_Q(X_m) \leq \langle R, T\rangle$.  Conversely, we have that $TR^m =M_{u_m}$ and $T R^{m+1}=M_{v_m}$.  Hence, $H_Q(X_m) = \langle R, T\rangle.$  Consequently, we have that $M_2 H_Q(X_m) M_2^{-1} = \langle S, T \rangle$.
\end{proof}

\begin{lemm}\label{lemm:STUGeneration}
We have $ST^{-1} S^{-1} T = - U_+^6$ and $T S^{-1} T^{-1} S= -U_-^6$. 
\end{lemm}
\begin{proof}  Regard $S, T$ as elements of $\SL(2, R_Q)$.
 We have that $ T(1, 0)^\mathsf{T} =  S (0, 1)^\mathsf{T} $ and $S^{-1}(1, 0)^\mathsf{T} = T^{-1}(0, -1)^\mathsf{T} $.  Consequently, \[-ST^{-1} S^{-1} T = \begin{pmatrix}  1 & b \\ 0 & 1 \end{pmatrix}.\]  Now computing $-ST^{-1} S^{-1} T (0, 1)^\mathsf{T}$ shows that $b=6$, giving the first equality.  Since $S = S^\mathsf{T}$, $T = T^\mathsf{T}$ and $U_- = U_+^\mathsf{T}$ hold, the other equality also follows.
\end{proof}

For a prime $p$ and an integer $e \geq 1$, define the mapping \[\pi_{p,e}: \PP^1(R_{p^e}) \rightarrow \PP^1(R_p), \quad [u: v] \mapsto [u \bmod p: v \bmod p]\] and the subgroups $\U_{p,e} :=\langle U_+^{pa}, U_-^{pa} : a \in R_{p^e} \rangle \leq \PGL(2, R_{p^e}).$  Note that $R_{p^e}$ is a local ring, and its unique maximal ideal is $(p)$ (see~\cite[Chapter~8]{AM69}).  Thus $\ker(R_{p^e} \rightarrow R_p) = (p)$.

\begin{lemm}\label{lemm:FiberTransivitiy}
We have that $\U_{p,e} x = \pi_{p,e}^{-1} \pi_{p,e} (x)$ for every $x \in \PP^1(R_{p^e})$. 
\end{lemm}
\begin{proof}  Let $x:=[u:v]$. Since we have $\pi_{p,e}\left(U_+^{pa}[u:v]\right) = \pi_{p,e}\left([u+vpa:v]\right) = [u  \bmod p:v  \bmod p]$ and the analog for $U_-^{pa}$, it follows that $\U_{p,e} x \subset \pi_{p,e}^{-1} \pi_{p,e} (x)$.

Let $y:=[u':v'] \in  \pi_{p,e}^{-1} \pi_{p,e} (x)$.  Thus, we have $[u' \bmod p: v' \bmod p]= [u \bmod p: v \bmod p]$.  There are two cases: $v \in R_{p^e}^*$ and $v \notin R_{p^e}^*$.  First consider the case $v \in R_{p^e}^*$.  Since $R_{p^e}$ is a local ring, we have the image of $v$ is a unit under the canonical ring homomorphism and thus $v' \bmod p$ is in $R_p^*$ and $v' \in R_{p^e}^*$.  Hence, $[u' (v')^{-1} \bmod p : 1\bmod p] = [u v^{-1}  \bmod p : 1\bmod p]$.  Consequently, there exists $k \in \ZZ$  such that $u' (v')^{-1} = u v^{-1} + pk$, giving $U_+^{pk} [u v^{-1}: 1] = [u' (v')^{-1}: 1]$.

Finally, consider the case $v \notin R_{p^e}^*$.  Since $R_{p^e}$ is a local ring, we have that $v \in (p)$ and, consequently, $u \in  R_{p^e}^*$.  Analogous to the previous case, we have that $u' \in  R_{p^e}^*$ and hence $[1 \bmod p : v' (u')^{-1}\bmod p] = [1 \bmod p : v u^{-1}  \bmod p]$.  Analogous to the previous case, we have that there exists $\ell \in \ZZ$  such that $U_-^{p\ell} [1:v u^{-1}] = [1:v' (u')^{-1}]$.  Together, the two cases give $\U_{p,e} x \supset \pi_{p,e}^{-1} \pi_{p,e} (x)$.
\end{proof}

\noindent Lemma~\ref{lemm:FiberTransivitiy} can be strengthened when  $p \geq 5$.
\begin{lemm}\label{lemm:TransForpGreat5}
For every prime $p \geq 5$, we have that both $\langle U_+^6, U_-^6 \rangle$ and $\langle S, T \rangle$ act transitively on $\PP^1(R_{p^e})$. 
\end{lemm}

\begin{proof}  By Lemma~\ref{lemm:STUGeneration}, it suffices to show $\langle U_+^6, U_-^6 \rangle$ acts transitively on $\PP^1(R_{p^e})$.  Since $6 \in R_{p^e}^*$, we have that, for every $a \in \ZZ$, there exists $k \in \ZZ$ such that $6k \equiv a \pmod{p^e}$ and thus $\langle U_+^6, U_-^6 \rangle = \langle U_+^a, U_-^a : a \in R_{p^e} \rangle$.  Since $R_{p^e}$ is a local ring, an arbitrary element is either $[u: 1]$ or $[1:u]$ for some $u \in R_{p^e}$.  Thus, we have $U_+^{-u+v} [u:1]= [v:1]$ and $U_-^{-u+v} [1:u]= [1:v]$.  Let $U:=U_+ U_- U_+^{-1}$ and note that $U[0:1] = [-1:0]=[1:0]$. Consequently, we have that $U_-^{v}UU_+^{-u}[u:1]=[1:v]$.
 
\end{proof}

\begin{prop}\label{prop:XmTransitive}
For every integer $Q \geq 2$, we have that $\langle S, T \rangle$ acts transitively on $\PP^1(R_Q)$. 
\end{prop}
 \begin{proof} Write $G:=\PGL(2, R_Q)$ and $Q = 2^a 3^b \prod_{j=1}^k p_j^{e_j}$ where $a, b \in \NN \cup \{0\}$ and $p_j \geq 5$ are distinct primes such that $p_j | Q$ with order $e_j:= \nu_{p_j}(Q):=\max\{e \in \NN: p_j^e | Q\}$.  By the Chinese remainder theorem, we have the ring isomorphism \[R_Q \cong \ZZ /2^a \ZZ \times \ZZ /3^b \ZZ \times \prod_{j=1}^k \ZZ / p_j^{e_j}\ZZ = R_{2^a} \times  R_{3^b} \times  \prod_{j=1}^k R_{p_j^{e_j}}\] and the bijection (\cite[Remark~8.31]{Ste07}) \[\PP^1(\R_Q) \cong \PP^1(R_{2^a}) \times  \PP^1(R_{3^b}) \times  \prod_{j=1}^k \PP^1(R_{p_j^{e_j}})\] where, throughout this proof, the factors corresponding to $a=0$ or $b=0$ are omitted.  Let $[u:v] \in \PP^1(R_Q)$ and $A \in G$ and define $[u:v]_p:= [u \mod p^e: v \mod p^e]$ and $A(p):= A \mod p^e$ where $p | Q$ is a prime and $e = \nu_p(Q)$.  Explicitly, the bijection gives \[[u:v] \mapsto \left([u:v]_{2}, [u:v]_{3}, [u:v]_{p_1}, \cdots,  [u:v]_{p_k}\right).\]  Note that the $A$-action is componentwise \[A [u:v] \mapsto \left(A(2) [u:v]_{2}, A(3) [u:v]_{3}, A(p_1) [u:v]_{p_1}, \cdots,  A(p_k )[u:v]_{p_k}\right).\]  Let $x:=(x_2, x_3, x_{p_1}, \cdots, x_{p_k}), y:= (y_2, y_3, y_{p_1}, \cdots, y _{p_k}) \in \PP^1(\R_Q)$.  We construct an element of $\langle S, T \rangle$ that takes $x$ to $y$ using the following steps.\\
 
\noindent {\bf Step~1:  the base modulo $2$ and $3$.}  Let us first assume that $a \geq 1$ and $b \geq 1$.  Note that $\{S^n [0:1]: n = 0, 1, 2\} =\PP^1(R_2)$  and that $\{[0:1], S[0:1], T[0:1], ST[0:1] \} = \PP^1(R_3)$.  Moreover, we have that $S^3 \equiv T \equiv  1_G \pmod 2$ and $S^2 \equiv 1_G \pmod 3$.  (Recall we are taking the projective image.)

Thus, for every $r \in \ZZ$, we have $\pi_{3,b}(S^{4r}(3) x_3) = \pi_{3,b}(x_3)$.  There exists $m \in \{0, 1, 2\}$ with $g_2:= S^{4m}$ such that $\pi_{2,a}(g_2(2) x_2) = \pi_{2,a}(y_2)$. Now there exists a word $g$ in the letters $S$ and $T$ such that $\pi_{3,b}(g(3) x_3) = \pi_{3,b}(y_3)$.  In $g$, replace every occurrence of $S$ with $S^3$ and $S^{-1}$ with $S^{-3}$ and call it $g_3$.  Thus, $g_3 \equiv 1_G \pmod 2$ and $g_3 \equiv g \pmod 3$.  Define $x^{(1)}:= g_3 g_2 x$.  Then we have that $x_2^{(1)}, y_2 \in \pi_{2,a}^{-1} \pi_{2,a} (y_2)$ and $x_3^{(1)}, y_3 \in \pi_{3,b}^{-1} \pi_{3,b} (y_3)$.  Note that, if $a=0$ or $b=0$, the corresponding part of Step~1 is omitted.  If both $a=0$ and $b=0$, Step~1 is omitted. \\

\noindent {\bf Step~2:  correcting the $2^a$-fiber.}  If $a=0$, this step is omitted.  Otherwise, assume that $a \geq 1$.  By Lemma~\ref{lemm:FiberTransivitiy}, there exists $h_2 \in \U_{2,a}$ such that $h_2 x_2^{(1)} = y_2$.  Let $U_+^{2 r}$ with $r \in R_{2^a}$ be a factor of $h_2$.  By the Chinese remainder theorem, there exists an $s \in R_Q$ such that \[6s \equiv 2 r \pmod {2^a}, \quad 6s \equiv 0 \pmod {Q/2^a}.\]  The analogous holds for $U_-^{2r'}$, yielding an $s' \in R_Q$.  Replace the factor $U_+^{2 r}$ with $U_+^{6s}$ and the factor  $U_-^{2r'}$ with  $U_-^{6s'}$.  Doing thus for all factors yields $\widetilde{h}_2$.  Consequently, applying Lemma~\ref{lemm:STUGeneration}, we have that $\widetilde{h}_2(2) x_2^{(1)} = y_2$ and $\widetilde{h}_2(p) x_p^{(1)} = x_p^{(1)}$ for all primes $p | Q$ such that $p \neq 2$.  Now define $x^{(2)}:= \widetilde{h}_2 x^{(1)}$ and note that $x^{(2)}_2 = y_2$ and $x^{(2)}_p = x^{(1)}_p$ for all primes $p | Q$ such that $p \neq 2$. \\

\noindent {\bf Step~3:  correcting the $3^b$-fiber.}  If $b=0$, this step is omitted.  Otherwise, assume that $b \geq 1$.  Analogous to Step~2, we obtain $\widetilde{h}_3(3) x_3^{(2)} = y_3$ and $\widetilde{h}_3(p) x_p^{(2)} = x^{(2)}_p$ for all primes $p | Q$ such that $p \neq 3$. Now define $x^{(3)}:= \widetilde{h}_3 x^{(2)}$ and note that $x^{(3)}_2 = y_2$, $x^{(3)}_3 = y_3$, and $x^{(3)}_p = x^{(2)}_p$ for all primes $p | Q$ such that $p \neq 2, 3$. \\

 \noindent{\bf Step~4:  correcting each $p_j^{e_j}$-component.}  If $k=0$, this step is omitted.  Otherwise, let $j=1$ and $K:= p_j^{e_j}$.   By Lemmas~\ref{lemm:TransForpGreat5} and~\ref{lemm:STUGeneration}, there exists $h_{j+3} \in \langle U_+^6, U_-^6 \rangle$ such that $h_{j+3}(p_j) x_{p_j}^{(j+2)} = y_{p_j}$.  Let $U_+^{6 r}$ with $r \in R_K$ be a factor of $h_{j+3}$.  By the Chinese remainder theorem, there exists an $s \in R_Q$ such that \[6s \equiv  6r \pmod K, \quad 6s \equiv 0 \pmod{Q/K}.\]  Consequently, $U_+^{6 s}$ acts trivially on the $Q/K$-component.  Analogously for a  factor $U_-^{6 r'}$ of $h_{j+3}$, there exists an $s' \in R_Q$.  Replacing the factors as in Step~3 yields a $\widetilde{h}_{j+3}$ such that $\widetilde{h}_{j+3}(p_\ell) x_{p_\ell}^{(j+2)} = y_{p_\ell}$ for all $1 \leq \ell \leq j$, $\widetilde{h}_{j+3}(2) x_2^{(j+2)} = y_2$, $\widetilde{h}_{j+3}(3) x_3^{(j+2)} = y_3$, and $\widetilde{h}_{j+3}(p_\ell) x_{p_\ell}^{(j+2)} = x_{p_\ell}^{(j+2)}$  for all $j < \ell \leq k$. Now define $x^{(j+3)}:= \widetilde{h}_{j+3} x^{(j+2)}$.  Finally, by recursion of Step~4 for $j = 2, \cdots k$, we obtain the desired result.

\end{proof}

\begin{prop}\label{prop:PosHD}
Let $Q \geq 2$ be an integer and $F\subset \QQ$ be nonempty and finite with $Qr \in \ZZ$ for all $r \in F$.  Then $d_F(t) >0$ for all $t > 3Q$. 
\end{prop}

\begin{proof}
 Fix $t > 3Q$ and put $s:= t/Q >3$.  Choosing $m$ large enough, we have $\mathfrak{H}(X_m) < s$.  Let \[\alpha_m := \min_{w \in X_m} \inf_{y \in [0,1]} |M'_w(y)| >0.\]  Then, for small enough $s_m >0$, Lemma~\ref{lemm:pressure}~(ii) gives $\dim_H K(X_m) >0$.  The result follows by applying Lemma~\ref{lemm:STConjHQXm}, Proposition~\ref{prop:XmTransitive}, and Theorem~\ref{prop:transitivesteering}.
\end{proof}

\begin{rema*}
The strong approximation method from Section~\ref{subsec:FullHD} is not quite sufficient for this section.  For example, $H_2(X_m)$ is a proper subgroup of $\PGL(2, R_2)$.

\end{rema*}

\subsection{Scaling law and full Hausdorff dimension}   In this section, we prove two of the main results of this paper, Theorem~\ref{thm:MainBad}~\ref{it3:Scaling:thm:MainBad} and~\ref{it2:FullHD:thm:MainBad},  using transitive projective cocycle steering.  Arithmetic enters strongly again through Proposition~\ref{prop:XmTransitive}. We also use symbolic arguments involving connections and branching; related constructions appear in \cite[Section~7.1]{MM20}, \cite[Section~3.2]{MMPV22}, and \cite[Section~2]{EJM25}. We combine these ideas with a return-word construction to obtain alphabets satisfying the required height and dimension bounds while retaining transitivity of the projective cocycle.  (Note that our return word is slightly different from the usual definition in, for example,~\cite{Dur98}.)  First, we choose a suitable prefix-free alphabet $X_*:=\{11,22\}$, which is a subset of the alphabet introduced in Section~\ref{subsec:FullHD} (see also~\cite[Section~7.1]{MM20}).  Let \[ x_1:= [0; 1, 1, \cdots] = \frac {\sqrt 5 -1} 2, \quad x_2:= [0;2,2,\cdots] = \sqrt 2- 1, \quad s_*:=2 + x_1 + x_2 = \frac{1 + 2 \sqrt 2 + \sqrt 5}{2}\] and set \[\widetilde{\xi}_1 = (1)_{j \in \NN}, \quad \xi_1:=(1)_{j \in \ZZ}   \in \Sigma(X_*), \quad \widetilde{\xi}_2 = (2)_{j \in \NN}, \quad \xi_2:=(2)_{j \in \ZZ} \in \Sigma(X_*).\]   Also let \[\omega:= \widetilde{\xi}_1^\mathsf{T} 2 2 1 1 \widetilde{\xi}_2 \in \Sigma(X_*), \quad \psi:= \widetilde{\xi}_1^\mathsf{T} 22 \widetilde{\xi}_1 \in \Sigma(X_*)\] where the $0$-th position of both bi-infinite sequences is the first $2$ from the left.

\begin{lemm}\label{lemm:XstarHeightsDims}
We have that $\dim_H K(X_*) >0$, $\max K(X_*) = x_1$, $\min K(X_*)= x_2$, $\mathfrak{H}(X_*) =s_*$, $\mathfrak{m}(\omega) < s_*$, $\mathfrak{m}(\mathcal{R}\omega) < s_*$, and $\mathfrak{m}(\psi) < s_*$.  Moreover, the only elements $\theta$ of $\Sigma(X_*)$ such that $f(\theta) = s_*$ are $\theta=\widetilde{\xi}_1^\mathsf{T} \widetilde{\xi}_2$ and $\theta = \mR(\widetilde{\xi}_1^\mathsf{T} \widetilde{\xi}_2)$.
\end{lemm}

\begin{proof}  The first assertion follows in the analogous way to that in Proposition~\ref{prop:PosHD}.  The fixed points of $M_{11}$ and $M_{22}$ in $[0,1]$ are $x_1$ and $x_2$, respectively.  Since both maps are increasing contractions on $[0,1]$, we have $M_{11}([x_2, x_1]) \subset [x_2, x_1]$ and $M_{22}([x_2, x_1]) \subset [x_2, x_1]$.  Thus, $K(X_*) \subset[x_2, x_1]$.  As $x_1, x_2 \in K(X_*)$, the second and third assertions follow.

Now let $\theta\in \Sigma(X_*)$.  If $\theta_0=1$, then $f(\theta) < 3$.  Otherwise, $\theta_0 =2$.  Then either $\theta_{-1}=2$ or $\theta_1 = 2$.  In either case, we have that $f(\theta) = 2 + u + \frac 1 {2 + v} \leq 2 + x_1 +   \frac 1 {2 + x_2} =s_*$ for some $u, v \in [x_2, x_1]$.  Thus, for $\theta':= \widetilde{\xi}_1^\mathsf{T} \widetilde{\xi}_2 \in \Sigma(X_*)$ with the $0$-th position at first $2$ from the left, we have $f(\theta') = s_*$, yielding the fourth assertion.  Moreover, the only other element $\theta'' \in \Sigma(X_*)$ such that $f(\theta'') = s_*$ is $\theta'' = \mathcal{R} \theta'$.

To show the final three assertions, it suffices to exclude equality.  The orbit closures $\{\sigma^n \omega: n \in \ZZ\} \cup \{\xi_1, \xi_2\}$, $\{\sigma^n \mathcal{R} \omega: n \in \ZZ\} \cup \{\xi_1, \xi_2\}$, and $\{\sigma^n \psi: n \in \ZZ\} \cup \{\xi_1\}$ are $\sigma$-invariant compact subsets that are disjoint from $\{\theta', \theta''\}$ and, since $f$ is continuous, equality is excluded in all three assertions.
\end{proof}

\subsubsection{Connections and dimension}

Let $(F_r)$ be the Fibonacci sequence with $F_1 = 1= F_2$ and define $\epsilon(r):= F_r^{-2}$.

\begin{lemm}\label{lemm:RightTails}
Let $n \geq 2$ be an integer, $M >0$, and $\alpha, \beta, \gamma \in \Sigma(\{1,2\})$ be such that $\mathfrak{m}(\alpha) \leq M, \mathfrak{m}(\beta) \leq M$, $\alpha_j = \beta_j$ for all $0 \leq j \leq n$, and \begin{align}\label{eqn:lemm:RightTails}\gamma_j = \beta_j \quad (j \leq n), \quad \gamma_{n+m} = \zeta_m \quad (m \geq 1)
  \end{align} for $\zeta_m \in \{1,2\}$ with \[\min([0;\alpha_{n+1}, \cdots ], [0;\beta_{n+1}, \cdots]) \leq  [0; \zeta_1, \cdots] \leq \max([0;\alpha_{n+1}, \cdots ], [0;\beta_{n+1}, \cdots]).\]  Then \[f(\sigma^j \gamma) \leq M + \epsilon(\lfloor n/2 \rfloor+1) \] for all $j \leq n$.  The same result holds if $\beta$ is replaced by $\alpha$ in (\ref{eqn:lemm:RightTails}).
\end{lemm}

\begin{proof}
 By~\cite[Lemma~3.6]{MMPV22}, the desired result follows if $\alpha_j = \beta_j$ for all $j \leq 0$.  Otherwise, as in the proof of Lemma~\ref{lemm:transSWthres}, we have an additional error term $ \epsilon(\lfloor n/2 \rfloor+1)$, giving the desired result.
\end{proof}

Let $s>0$, let $Z \subset \Sigma(\{1,2\})$ be nonempty, compact, and $\sigma$-invariant, and define \[Z^+:=\{(\theta_0, \theta_1, \cdots): \theta \in Z\}, \quad \mathfrak{H}(Z) := \max_{\theta\in Z}f(\theta).\]   A {\em right connection from $Z$ to $\xi_i$ below $s$} is an element $\gamma \in \Sigma(\{1,2\})$ with $\mathfrak{m}(\gamma) <s$ for which some $\theta \in Z$ and integers $\ell < r$ satisfy \[\gamma_j = \theta_j \quad (j \leq \ell), \quad \gamma_j = i \quad (j \geq r).\]  Likewise,  a {\em left connection from $\xi_i$ to $Z$ below $s$} is an element $\gamma \in \Sigma(\{1,2\})$ with $\mathfrak{m}(\gamma) <s$ for which some $\theta \in Z$ and integers $\ell < r$ satisfy  \[ \gamma_j = i  \quad (j \leq \ell), \quad  \gamma_j = \theta_j \quad (j \geq r).\]  
 
Let $\widetilde{Z} \subset \{1, 2\}^{\NN \cup \{0\}}$. A word $w$ is {\em a prefix of $\widetilde{Z} $} if $w$ is a prefix of some element of $\widetilde{Z} $ and is a {\em branching prefix (of $\widetilde{Z}$)} if both $w1$ and $w2$ are prefixes of $\widetilde{Z} $.  Such a branching prefix has the {\em pair property} if every element of $\widetilde{Z} $ beginning with $w1$ begins with $w11$, and every element of $\widetilde{Z} $ beginning with $w2$ begins with $w22$. Finally, for $\emptyset \neq K \subset K(\{1,2\})$, define \[\mZ(K) := \{\theta \in \Sigma(\{1,2\}): [0;\theta_0, \theta_1, \cdots ] \in K\}, \quad \mZ^+(K):= \{(\theta_0, \theta_1, \cdots): \theta \in \mZ(K)\}.\]

\begin{lemm}\label{lemm:NoConnection}
Let $s> \max(\mathfrak{H}(Z), 3)$ and $r_0 := \min\{r \in \NN \backslash\{1\}: \epsilon(r) < s- \mathfrak{H}(Z)\}$.  If $Z$ has no right connection below $s$ to either $\xi_1$ or $\xi_2$, then every branching prefix $w$ of $Z^+$ with $|w|_{\NN} \geq 2r_0-1$ has the pair property.
\end{lemm}
\begin{proof}
Let $w:=a_0 \cdots a_n$ be a branching prefix of $Z^+$ where $n \geq 2 r_0 - 2$.  Then there are elements of $Z^+$ for which $w1$ is a prefix and elements for which $w2$ is a prefix.  Assume that the conclusion does not hold.  We have two cases.\\

\noindent{\bf Case:  $w12$ is a prefix.}  Let $\alpha, \beta \in Z$ such that $\alpha_0 \cdots \alpha_n = w = \beta_0 \cdots \beta_n$, $\alpha_{n+1} = 1$, $\alpha_{n+2} = 2$, and $\beta_{n+1}=2$.  Define $\gamma \in \Sigma(\{1,2\})$ by \[\gamma_j = \begin{cases}  \beta_j & \textrm{ if } j \leq n \\ 1 &\textrm{ if } j > n\end{cases}.\]  Since $[0;\beta_{n+1}, \cdots]<\frac12<x_1<\frac23<[0; \alpha_{n+1}, \cdots]$, we may apply Lemma~\ref{lemm:RightTails} to obtain $f(\sigma^j \gamma) \leq \mathfrak{H}(Z) + \epsilon(\lfloor n/2 \rfloor+1) \leq \mathfrak{H}(Z) + \epsilon(r_0) < s$ for all $j \leq n$.
  For $j >n$, we have $f(\sigma^j \gamma) < 3$. Thus, $\gamma$ is a right connection from $Z$ to $\xi_1$ below $s$, a contradiction.\\
  
  \noindent{\bf Case:  $w21$ is a prefix.}  Let $\alpha, \beta \in Z$ such that $\alpha_0 \cdots \alpha_n = w = \beta_0 \cdots \beta_n$, $\alpha_{n+1} = 2$, $\alpha_{n+2} = 1$, and $\beta_{n+1}=1$.  Define $\gamma \in \Sigma(\{1,2\})$ by \[\gamma_j = \begin{cases}  \alpha_j & \textrm{ if } j \leq n \\ 2 &\textrm{ if } j > n\end{cases}.\]  Analogous to the previous case, we obtain  $f(\sigma^j \gamma) < s$ for all $j \leq n$.  Now, since $x_2 <  [0; \alpha_{n+2}, \cdots]$, we have $f(\sigma^{n+1} \gamma) \leq f(\sigma^{n+1} \alpha) <s$.  For $j \geq n+2$, we have $f(\sigma^j \gamma) \leq 2 + x_2 + 1/2 < 3$. Thus, $\gamma$ is a right connection from $Z$ to $\xi_2$ below $s$, a contradiction.   \end{proof}
  
Let $w:=(a_1, \cdots, a_n)$ be a word over $\NN$.  Integrating the monotonic function from (\ref{eqn:DerivativeCocycle}) over $[0,1]$ yields $|\mathcal{I}(w)| = (q_n(q_n+q_{n-1}))^{-1}$, and setting \begin{align}\label{eqn:DefnbwContFrac}
 b_w:= \inf_{y \in [0,1]} |M'_w(y)| \end{align} and using elementary properties of continued fractions yields $b_w \leq |\mathcal{I}(w)| \leq 2 b_w$.
   
\begin{lemm}\label{lemm:PairPropLowDim} Let $K \subset K(\{1,2\})$ be nonempty and compact.  If there exists $N \in \NN$ for which every branching prefix $w$ of $\mZ^+(K)$ satisfying $\|w\|_{\NN} \geq N$ has the pair property, then $\dim_H(K) \leq \dim_H K(X_*)$.
 \end{lemm}
\begin{proof}  Let $n \in \NN$, fix $u > \dim_H K(X_*)$, choose $v$ with $\dim_H K(X_*) < v <u$, and let $w \in X_*^n$.  Since \[\sup_{y \in [0,1]} |M'_{\widetilde{w}}(y)|  \leq 4^{-\frac {|\widetilde{w}|_{\NN}}2}\] for $\widetilde{w} = 11$ and $\widetilde{w}=22$, the chain rule gives the same for $\widetilde{w} = w$.  Integration over $[0,1]$ yields $|\mathcal{I}(w)| \leq 4^{-n}$.   Consequently, we have that \begin{align}\label{eqn:lemm:PairPropLowDim:pressure}
 \sum_{w \in X_*^n} |\mathcal{I}(w)|^u \leq 2^v 4^{-n(u-v)}\sum_{w \in X_*^n} b_w^v \rightarrow 0 \end{align} as $n \rightarrow \infty$ by Lemma~\ref{lemm:pressure}~\ref{it2:lemm:pressure}.

Let $\mZ:= \mZ^+(K)$, $\mP_N$ denote the set of prefixes $w$ of $\mZ$ with $|w|_{\NN}=N$, and $\B$ denote the set of branching prefixes of $\mZ$.  Fix $w \in \mP_N$ and let $\mZ_w := \{\theta \in\mZ : \theta \textrm { begins with } w\}$.  For $m \in \NN$ and $\theta \in \mZ_w$, define $\Fp_m(\theta) := \theta_0 \theta_1 \cdots \theta_{m-1}$, $\mZ_{w,\infty} :=\{\theta \in \mZ_w: \Fp_m(\theta) \in \B \textrm{ for infinitely many } m\}$, and $K_{w, \infty}:= \{[0; \theta_0, \theta_1,  \cdots] : \theta \in \mZ_{w, \infty}\}$.  Let $\mH^u(\cdot)$ denote the $u$-dimensional Hausdorff measure.  We now show that $\mH^u(K_{w, \infty})=0$.  There are two cases.\\

\noindent{\bf Case:  $\mZ_{w,\infty} \neq \emptyset.$}  Recursively, define $t_1(\theta) := \min \{ m \geq N: \Fp_m(\theta) \in \B\}$ and $t_{j+1}(\theta) := \min \{ m \geq t_j(\theta)+2: \Fp_m(\theta) \in \B\}$ for all $j \in \NN$.  Put $\alpha_j(\theta):= \theta_{t_j(\theta)} \in \{1, 2\}$.  Then the pair property gives $\theta_{t_j(\theta)}\theta_{t_j(\theta)+1} = \alpha_j(\theta)\alpha_j(\theta) \in X_*$.  Now define \[\Fb_n(\theta):= \alpha_1(\theta)\alpha_1(\theta) \cdots \alpha_n(\theta)\alpha_n(\theta)\in X_*^n,\quad \mW_n(w):=\{\Fp_{t_n(\theta)+2}(\theta): \theta \in \mZ_{w,\infty}\},\] and \[\Tr_{w,n}: \mW_n(w) \rightarrow X_*^n, \quad \Fp_{t_n(\theta)+2}(\theta) \mapsto  \Fb_n(\theta).\]  

We claim the mappings $\Tr_{w,n}$ are injections. We now prove this claim.  Let $\Tr_{w,n}(\widehat{w}) = \Tr_{w,n}(\widecheck{w})$.  Then $w$ is a prefix of both $\widehat{w}$ and $\widecheck{w}$.  Starting with $w$, we reconstruct the successive digits of $\widehat{w}$ and $\widecheck{w}$.  For a prefix $\nu$ of $\mZ$ that is not branching, we have either $\nu 1$ or $\nu 2$ is a prefix of $\mZ$ but not both, and, consequently, the digit succeeding $\nu$ is uniquely determined.  For a branching prefix $\nu$ of $\mZ$, the next two digits are the same for both and determined by $\Fb_n(\theta)$.  This proves the claim.

Consequently, one obtains $\Fb_n(\theta)$ by removing digits from $\Fp_{t_n(\theta)+2}(\theta)$.  Let $q(\cdot)$ denote the continuant, and, by Euler’s property on continuants (\cite[Appendix~2]{CF89}) and elementary properties of continued fractions, we have $q(\Fb_n(\theta)) \leq q(\Fp_{t_n(\theta)+2}(\theta))$ and, thus, \begin{align}\label{eqn:lemm:PairPropLowDim:cover}
|\mathcal{I}(\Fp_{t_n(\theta)+2}(\theta))|  \leq  q(\Fp_{t_n(\theta)+2}(\theta))^{-2} \leq q(\Fb_n(\theta))^{-2} \leq  2 |\mathcal{I}(\Fb_n(\theta))| . \end{align}  Then, by construction, we have \[K_{w, \infty} \subset \bigcup_{w' \in \mW_n(w)} I(w')\] for all $n \geq 1$.  By~(\ref{eqn:lemm:PairPropLowDim:pressure}) and~(\ref{eqn:lemm:PairPropLowDim:cover}), we have that $\mH^u(K_{w, \infty})=0$.\\

\noindent{\bf Case:  $\mZ_{w,\infty} = \emptyset.$}  Then $[K_{w, \infty}=\emptyset$ and thus, $\mH^u(K_{w, \infty})=0$.\\

Now let $\mZ_{\fin}$ be the set of elements in $\mZ$ containing only finitely many branching prefixes and define $K_{\fin}:= \{[0;\theta_0, \theta_1, \cdots ] : \theta \in \mZ_{\fin}\}$.  Since the digits following a prefix that is not branching is uniquely determined, it follows that, after some long enough prefix, all subsequent digits are uniquely determined, giving that $K_{\fin}$ is countable.  Consequently, we have that \[K = K_{\fin} \cup \bigcup_{w \in \mP_N} K_{w, \infty}.\]  Note that $\mP_N$ has at most $2^N$ elements.  Hence, $\mH^u(K)=0$, yielding the desired result.\\
 
\end{proof}

% \subset \{1,2\}^+
\begin{lemm}\label{lemm:ExistConnections}
Let $X$ be a finite, prefix-free binary alphabet such that  $s> \max(\mathfrak{H}(X), 3)$ and $\dim_H K(X) > \dim_H K(X_*)$.  Then there exists a right connection below $s$ from $\Sigma(X)$ to some $\xi_i$ and there exists a left connection below $s$ from some $\xi_j$ to $\Sigma(X)$ where $i, j \in \{1, 2\}$.
\end{lemm}
\begin{proof}
 The first assertion follows from Lemmas~\ref{lemm:NoConnection},~\ref{lemm:PairPropLowDim}, and~\ref{lemm:AlphabetShiftSameCantorDim}.  Replacing $\Sigma(X)$ by $\mR\Sigma(X)$ in the argument for the first assertion, yields the second.
\end{proof}

\subsubsection{First return words and an associated alphabet}  Let $R \in \NN$ and $\Omega \subset \Sigma(\{1,2\})$ be nonempty, compact, and $\sigma$-invariant.  Define \begin{align*}
\eL_R(\Omega)&:= \{\theta_k \theta_{k+1} \cdots \theta_{k+R-1}: \theta \in \Omega, k \in \ZZ\}, \\ \Omega_R&:= \{\theta \in \Sigma(\{1,2\}): \theta_k \theta_{k+1} \cdots \theta_{k+R} \in \eL_{R+1}(\Omega) \textrm { for all } k \in \ZZ\}.\end{align*}  Note that $\eL_R(\Omega)$ is the set of all $R$-blocks that occur in at least one element of $\Omega$ (see~\cite[Definition~1.3.1]{LM21}).   A word $w:=a_1 \cdots a_n$ over $\{1, 2\}$ with length $n \geq R+1$ is {\em locally admissible (with respect to $R$ and $\Omega$)} if $a_k a_{k+1} \cdots a_{k+R} \in \eL_{R+1}(\Omega)$ for every $1 \leq k \leq n -R$.  We also refer to every word over $\{1, 2\}$ with length $n \leq R$ as {\em locally admissible (with respect to $R$ and $\Omega$)}.  With this, note that any nonempty subword of a locally admissible word with respect to $R$ and $\Omega$ is also locally admissible with respect to $R$ and $\Omega$.  Note that $\eL_R(\Omega) \subset \{1, 2\}^R$ is finite and that $\Omega_R \supset \Omega$ and is compact and $\sigma$-invariant.  Also note that, for $n \geq R+2$, local admissibility does not necessarily imply that $w$ occurs in an element of $\Omega$.

\begin{lemm}\label{lemm:LocalAdmissR} Let $s>0$.  If $\mathfrak{H}(\Omega) <s$, then there exists an even integer $R\geq 2$ such that $\Omega \subset \Omega_R$ and $\mathfrak{H}(\Omega_R) <s$.
\end{lemm}

\begin{proof}
 Choose $R$ such that $\mathfrak{H}(\Omega) +2 \epsilon_{R/2}<s$.  For $\theta \in \Omega_R$, we have that \[\theta_{-R/2} \theta_{-R/2+1} \cdots \theta_0 \theta_1 \cdots \theta_{R/2} = \eta_k \eta_{k+1} \cdots \eta_{k+R}\] for $\eta \in \Omega$ and thus $f(\theta) \leq f(\sigma^{k+R/2} \eta) + 2\epsilon_{R/2} <s$.
\end{proof}

For a word $w':=a'_1 \cdots a'_m$ over $\NN$ with $m \geq R$, define $\Fs_R(w') := a'_{m - R +1} \cdots a'_m$.   Now fix a word $\tau$ over $\{1, 2\}$ such that $|\tau|_{\NN} = R$.  A {\em return word (to $\tau$)} is a nonempty word $v:= v_1 \cdots v_\ell$ such that $\tau v$ is locally admissible with respect to $R$ and $\Omega$ and $\Fs_R(\tau v) = \tau$ and is, furthermore, a {\em first return word (to $\tau$)} if $\Fs_R(\tau v_1 \cdots v_k) \neq \tau$ for all $1 \leq k < \ell$.  We also refer to $\tau$ as the {\em reference word}.

\begin{lemm}\label{lemm:ReturnWordCat} Fix a reference word $\tau$ with $|\tau|_{\NN}=R$ and let $u, v$ be return words.  Then $uv$ is a return word.
 \end{lemm}
 
 \begin{proof}
 It suffices to show that every $(R+1)$-block in $\tau u v$ is in $\eL_{R+1}(\Omega)$.  Since $\Fs_R(\tau u) = \tau$ and $|\tau v|_{\NN}= R+ |v|_{\NN}$, we have that every $(R+1)$-block in $\tau u v$ is either in $\tau u$ or $\tau v$ and thus locally admissible.
\end{proof}

\begin{lemm}\label{lemm:FirstReturnWords}  Fix a reference word $\tau$ with $|\tau|_{\NN}=R$.  Then every return word is a concatenation of first return words, every finite alphabet $\F$ consisting of first return words is prefix-free, and $\Sigma(\F) \subset \Omega_R$.
 
\end{lemm}

\begin{proof}
 Let $v:=b_1 \cdots b_\ell $ be a return word and let $t_1 < t_2 < \cdots < t_k = \ell$ be the elements of $\mT:=\{t \in \{1, \cdots, \ell \}: \Fs_R(\tau b_1 b_2 \cdots b_t) = \tau\}$.  Set $t_0=0$ and, for $1 \leq j \leq k$, define $v_j:= b_{t_{j-1}+1} b_{t_{j-1}+2} \cdots b_{t_j}$.  Since the $R$ digits preceding $v_j$ in $\tau v$ is $\tau$, we have that $\tau v_j$ is a subword of $\tau v$ and thus is locally admissible.  Moreover, $\tau v_j$ ends in $\tau$.  Therefore the collection of $v_j$ are return words, and, by the definition of $\mT$, are first return words.
 
 Now it follows from the definition of first return word that a first return word cannot have a return word as a proper prefix.  Consequently, any finite alphabet consisting of first return words must be prefix-free.  Finally, since every $(R+1)$-block $\gamma$ of an element $\theta \in \Sigma(\F)$ must lie in a finite concatenation of the elements of $\F$, Lemma~\ref{lemm:ReturnWordCat} shows that $\gamma$ is locally admissible.  The analogous holds for any $\sigma^{k} \theta$, giving that $\Sigma(\F) \subset \Omega_R$.
  
 \end{proof}
 
 \begin{lemm}\label{lemm:ConstXhatFromX}  Let $X$ be a finite, prefix-free binary alphabet and $\F$ be a finite collection of first return words to $\tau$.  Suppose there exist $\lambda, \mu \in \{1,2\}^*$  such that for every $w \in X^+$ we have that \begin{align}\label{eqn1:lemm:ConstXhatFromX}
 \tau \lambda w \mu \textrm{ is locally admissible for } R \textrm{ and } \Omega, \quad \Fs_R(\tau \lambda w \mu) =\tau.\end{align}  Then, for every $0 < u < \dim_H K(X)$, there exists a finite, prefix-free alphabet $\widehat{X}$ consisting of first return words to $\tau$ such that \begin{align*}\label{eqn2:lemm:ConstXhatFromX}
\F \subset \widehat{X}, \quad \dim_HK(\widehat{X}) >u, \quad \Sigma(\widehat{X}) \subset \Omega_R.  \end{align*}
 
\end{lemm}
\begin{rema*}
Note that $\lambda$, $\mu$, or both may be the empty and that $\F$ may be empty.

\end{rema*}
 \begin{proof}
 For $N \in \NN$, define $\mV_N:= \{\lambda w \mu: w \in X^N\}$.  Since $X^N$ is prefix-free, $\mV_N$ is prefix-free.  Recall the definition of $Z_N^-(X, u)$ from (\ref{eqn:PartitionSumsPressure}) and the definition of $b_w$ from (\ref{eqn:DefnbwContFrac}).  Let $b_\emptyset =1$.  By Lemma~\ref{lemm:pressure}~\ref{it1:lemm:pressure}, there exists an $N_0 \in \NN$ such that, for all $N \geq N_0$, we have \[\sum_{\nu \in \mV_N}b^u_\nu \geq b^u_\lambda b^u_\mu Z_N^-(X, u) > 1.\]  Fix such an $N$.  Applying Lemma~\ref{lemm:pressure}~\ref{it2:lemm:pressure} yields $u < \dim_H K(\mV_N)$.
 
 By (\ref{eqn1:lemm:ConstXhatFromX}), we have that every $\nu \in \mV_N$ is a return word to $\tau$, and, by Lemma~\ref{lemm:FirstReturnWords}, $\nu$ is a concatenation of first return words.  Let $\F(\nu)$ be the collection of these first return words and let \[\F_N := \bigcup_{\nu \in \mV_N} \F(\nu), \quad \widehat{X}:= \F_N \cup \F.\]  Note that $\Sigma^+(\mV_N) \subset \Sigma^+(\widehat{X})$ and hence $K(\mV_N) \subset K(\widehat{X})$.  Applying Lemma~\ref{lemm:FirstReturnWords} again yields the desired result.

\end{proof}

\begin{prop}\label{prop:XhatExistGoodProps}  Let $s \geq s_*$, $X$ be a finite, prefix-free binary alphabet with $\mathfrak{H}(X) < s$ and $\dim_H K(X) > \dim_H K(X_*)$.  For every $0 < u < \dim_H K(X)$, there exists a finite, prefix-free alphabet $\widehat{X}$ such that $\mathfrak{H}(\widehat{X}) < s$, $\dim_H K(\widehat{X})>u$, and $H_Q(\widehat{X})$ acts transitively on $\PP^1(R_Q)$ for every integer $Q \geq 2$.

\end{prop}
 
 \begin{rema*}
Note that $\widehat{X}$ does not depend on $Q$.

\end{rema*}
\begin{proof}
 Let $Z:= \Sigma(X)$.  Let $\eta^-$ be the left connection below $s$ from $\xi_i$ to $Z$ and $\eta^+$ be the right connection below $s$ from $Z$ to $\xi_j$ given by Lemma~\ref{lemm:ExistConnections}.  Here $i,j \in \{1,2\}$ and need not be equal.
 
 For $\theta \in \Sigma$, define $\mO(\theta):= \{\sigma^k \theta: k\in \ZZ\}$ and let $\overline{\mO(\theta)}$ denote its closure.  Let \[\Omega := \Sigma(X) \cup \overline{\mO(\eta^-)} \cup \overline{\mO(\eta^+)}\cup \overline{\mO(\omega)} \cup \overline{\mO(\mR\omega)}\cup \overline{\mO(\psi)}.\]  Then $\Omega$ is a finite union of compact and $\sigma$-invariant subsets and is thus compact and $\sigma$-invariant.  Lemma~\ref{lemm:XstarHeightsDims} and continuity of $f$ give that $\mathfrak{H}(\Omega) <s$.   Note the strict inequality holds for $s= s_*$ too.
 
 Let $R$ be the even integer from Lemma~\ref{lemm:LocalAdmissR}.  Local admissibility in the proof of this proposition is with respect to $R$ and $\Omega$.  Set the reference word to be $\tau:= 1^R$.  To allow us to construct return words to $\tau$, need to ensure that our connections are to and from $\xi_1$.  We do this as follows.
 
If $i=1$, then take $\widetilde{\eta}^-:= \eta^-$.  Otherwise, $i=2$.  Choose $\ell \in \ZZ$ such that $\eta^-_k = 2$ for all $k \leq \ell$ and define \[\widetilde{\eta}^-_k := \begin{cases}  \omega_{k+3 +R} & \text{ if }  k \leq 0 \\ \eta^-_{k+ \ell -R} & \text{ if } k \geq 1\end{cases}.\]  Note that $\widetilde{\eta}^-_{[-R+1, R]} = 2^{2R}$.  If an $(R+1)$-block contains both positions $0$ and $1$, then it is $2^{R+1}$.  Otherwise, it is either contained in $\omega$ or $\eta^-$.  Consequently, every $(R+1)$-block is in $\eL_{R+1}(\Omega)$ and thus $\widetilde{\eta}^- \in \Omega_R$.  Consequently, $\widetilde{\eta}^-$ is a left connection below $s$ from $\xi_1$ to $Z$.  We have an analogous construction for $\widetilde{\eta}^+$ with $\omega$ replaced by $\mR \omega$.

Now there exists some $m \in \NN$ such that $\widetilde{\eta}^-_{[m, \infty]} = \alpha_1 \alpha_2 \cdots$ where $\alpha_j \in X$.  Choose an occurrence of $\tau$ in $\widetilde{\eta}^-_{[-\infty, 0]}$ and form the subword $w^*$ of $\widetilde{\eta}^-$ beginning with this occurrence of $\tau$ and ending in the subword $\alpha:= \alpha_1 \cdots \alpha_R$ of $\widetilde{\eta}^-_{[m, \infty]}$ with the first position of $\alpha$ equal to the $m$-th position of $\widetilde{\eta}^-$:  $w^* = \tau \lambda \alpha$.  Thus we have that $\alpha \in X^+$ and $|\alpha|_{\NN} \geq R$.  Since $w^*$ is a subword of $\widetilde{\eta}^-$, it is locally admissible.  Analogously, using $\widetilde{\eta}^+$, we have $\beta \mu \tau$ and that it is locally admissible. Moreover, $\beta \in X^+$ and $|\beta|_{\NN} \geq R$.

We claim that, for every word $w \in X^+$, $\tau \lambda \alpha w \beta \mu \tau$ is locally admissible.  We now prove this claim.  Note that $\alpha w \beta$ is a concatenation of words in $X$ and thus occurs as a subword of an element of $\Sigma(X)$ and hence is locally admissible.  Now any $(R+1)$-block in $\tau \lambda \alpha w \beta \mu \tau$ must occur in $\tau \lambda \alpha$, $\alpha w \beta$, or $\beta \mu \tau$, proving the claim.

Now both $1$ and $22 \tau$ are first return words.  Let $\F:= \{1, 22 \tau\}$ and $\widehat{X}$ be the alphabet given by Lemma~\ref{lemm:ConstXhatFromX}.  Since Lemmas~\ref{lemm:ConstXhatFromX} and~\ref{lemm:LocalAdmissR} hold, we have that $\mathfrak{H}(\widehat{X}) \leq \mathfrak{H}(\Omega_R) <s$.  Note that $S, T \in H_Q(\F) \subset H_Q(\widehat{X})$ and, thus, by Proposition~\ref{prop:XmTransitive}, we have $H_Q(\widehat{X})$ acts transitively on $\PP^1(R_Q)$ for every integer $Q \geq 2$.
\end{proof}

\subsubsection{Transitive Cantor-set dimension function}  Let $Q \geq 2$ be an integer.  Using (\ref{eqnPerronConv}), let $D(s) = \dim_H(k^{-1}(-\infty, s))$ (see~\cite[Theorem~1(i)]{Mor18}).  By~\cite[Lemma~2]{Mor18}, we have that \[D(s) = \sup \left\{\dim_H K(X): X \textrm{ is a finite, prefix-free alphabet}, \mathfrak{H}(X) < s\right\}.\]   Now define \[D^{\tra}_Q(s) := \sup \left\{\dim_H K(X): \begin{array}{l} X \textrm{ is a finite, prefix-free alphabet}, \mathfrak{H}(X) < s, \\ H_Q(X) \textrm{ acts transitively on }\PP^1(R_Q)\end{array}\right\}.\]  For both definitions, if no finite, prefix-free alphabet $X$ exists, we take the convention that the supremum is $0$.

\begin{lemm}\label{lemm:TranDimLowerBnd}
Let $0< \varepsilon < \dim_H(K(X_*))$. Then there exists an integer $n \geq 2$ such that, for the finite, prefix-free alphabet $Y_* := X_*^n  \backslash \{(11)^n,(22)^n\}$, we have \[\dim_HK(Y_*) > \dim_H(K(X_*)) - \varepsilon, \quad \mathfrak{H}(Y_*) < s_*,\] and $H_Q(Y_*)$ acts transitively on $\PP^1(R_Q)$.  Moreover, we have $D^{\tra}_Q(s_*) \geq \dim_HK(X_*)$.
\end{lemm}

\begin{proof}  By Lemma~\ref{lemm:XstarHeightsDims}, we have $\dim_H(K(X_*))>0$.  Then Lemma~\ref{lemm:HighXRemoveExtremes} and Proposition~\ref{prop:XmTransitive} give the first and the third assertions.  Now every element of $Y_*$ contains both $11$ and $22$ as subwords.  Thus, Lemma~\ref{lemm:XstarHeightsDims} gives the second assertion.  Since $X_*$ is finite and prefix-free, so is $Y_*$.  Letting $\varepsilon \rightarrow 0$ yields the fourth assertion. \end{proof}

\begin{prop}\label{prop:ScalingLaw}  For every integer $Q \geq 2$, we have that \[D^{\tra}_Q(s) = D(s)\] for all $s_* \leq s \leq \sqrt{12}$.  Consequently, for every nonempty, finite $F \subset \QQ$ with $Qr \in \ZZ$ for all $r \in F$, we have \[d_F(t) \geq d_{\cl}(t/Q)\] for all  $t \geq Qs_*$.
 
\end{prop}

\begin{proof}  We have that $\mathfrak{m}(\theta) \leq \sqrt{12}$ if and only if $\theta \in \Sigma(\{1,2\})$~\cite[Page~35]{CF89}.  Thus every finite, prefix-free alphabet $X$ with $\mathfrak{H}(X) < s$ is binary.  By monotonicity and Lemma~\ref{lemm:TranDimLowerBnd}, we have $D^{\tra}_Q(s)\geq D^{\tra}_Q(s_*) \geq  \dim_H K(X_*).$   If $\dim_H(K(X_*)) \geq \dim_H K(X)$, then we have \begin{align}\label{eqn1:prop:ScalingLaw} D^{\tra}_Q(s)\geq\dim_H K(X).
  \end{align}

Otherwise, we have $\dim_H(K(X_*)) < \dim_H K(X)$.  Then, for every $0 < u < \dim_H K(X)$, Proposition~\ref{prop:XhatExistGoodProps} gives a suitable finite, prefix-free alphabet $\widehat{X}$ such that $\dim_H(K(\widehat{X}))>u$ and, consequently, $D^{\tra}_Q(s)\geq u$.  Letting $u \rightarrow \dim_H K(X)$ yields (\ref{eqn1:prop:ScalingLaw}).  This shows the first assertion.

Theorem~\ref{prop:transitivesteering} and~\cite[Theorem~1(i)]{Mor18} imply that $d_F(Qs) \geq \min(1, 2 D^{\tra}_Q(s)) = \min(1, 2 D(s)) = d_{\cl}(s)$ for all $s_* \leq s\leq \sqrt{12}$.   Note that~\cite[Theorem~1(iii)]{Mor18} gives $d_{\cl}(\sqrt{12})=1$ and thus, $d_F(Q\sqrt{12})=1$.  For $s > \sqrt{12}$, monotonicity gives $d_F(Qs) \geq 1 = d_{\cl}(s)$.  Setting $t := Qs$ yields the desired result.
\end{proof}

\begin{coro}\label{cor:FullHD} For every integer $Q \geq 2$ and every nonempty, finite $F \subset \QQ$ with $Qr \in \ZZ$ for all $r \in F$, we have $d_F(t) = 1$ for all $t \geq Q \Ft_1$.
\end{coro}
\begin{proof}
 Since $\Ft_1 \approx 3.33$~\cite[Theorem~2.1]{MMPV22} is greater than $s_*$ and $d_{\cl}(\cdot)$ is continuous~\cite[Theorem~1(i)]{Mor18}, the result follows from Proposition~\ref{prop:ScalingLaw}. 
\end{proof}

\section{A simultaneous Hall ray for rational translations}\label{sec:HallRay}   In this section, we prove one of the main results of this paper Theorem~\ref{thm:MainBad}~\ref{it1:HallRay:thm:MainBad}.  Our proof relies on Hall's theorem \eqref{eqn:HallThm}. The use of increasingly long windows joined by congruence-controlling words also appears in \cite[Section~2]{EM26}. Our setting, however, requires equality of Lagrange constants simultaneously under all rational translations with a fixed denominator. To achieve this, we develop Hall projective cocycle steering, introduced in Lemma~\ref{lemm:HallSteering}.  Arithmetic enters via Lemma~\ref{lemm:LagrangeUnderQPush} and Proposition~\ref{prop:XmTransitive}.  Let $s_\#:= 4 \sqrt 2$.  We have the finite, prefix-free alphabet $X:= \{1, 2, 3, 4\}$, and, for $M \in \NN$, let $C_M := K(\{1, \cdots, M\})$.  In the continued fraction notation, an overline to denotes the standard convention of infinite repeating blocks (see~\cite[Page~48]{Kh}).

\begin{lemm}\label{lemm:HallRayBasics}  We have $\min C_M = [0; \overline{M, 1}] = (\sqrt{M^2+4M} \ -M)/(2M)$, $\max C_M = [0; \overline{1, M}] = (\sqrt{M^2+4M} \ -M)/2$, $\mathfrak{H}(X) = s_\#$, and the only element $\theta \in \Sigma(X)$ such that $f(\theta) = s_\#$ is $\theta = \cdots 1 4 1 4 1 4 \cdots$ with a $4$ in the $0$-th position.  Note that $\min C_4 =  (\sqrt 2 - 1)/2$, $\max C_4 =2 (\sqrt 2 -1)$, and $s_\# = 4 + 2 \max C_4$.
 
\end{lemm}

\begin{proof}  By~\cite[Chapter~1, Lemma~1]{CF89} and directly computing, we have the first two assertions.  The proofs of the other assertions are analogous to that of Lemma~\ref{lemm:XstarHeightsDims}.
\end{proof}

We also have Hall's theorem (\cite{Hal}, \cite[Page~146]{Mor18}) \begin{align}\label{eqn:HallThm} C_4 +C_4 =[\sqrt 2 - 1, 4(\sqrt 2 - 1)]=:\eJ.
  \end{align}
  
\begin{lemm}\label{lemm:AttainHeightWithBigCenter}  For every $s \geq s_\#$, there exists an integer $a_0 \geq 5$ and $a, b \in \Sigma^+(X)$ such that \[\theta^{(s)}:= \cdots b_3 b_2 b_1 a_0 a_1 a_2 a_3 \cdots \in \Sigma, \quad  \theta^{(s)}_0 = a_0\] satisfies $f(\theta^{(s)}) = \mathfrak{m}(\theta^{(s)}) = s$.
 
\end{lemm}

%, and $\mathfrak{r}:=1, 4, 1, 4, \cdots \in \Sigma^+(X)$
\begin{proof}
 Let $c := \sqrt 2 - 1$, $r:=[0; \overline{1, 4}]$, and $\chi:=1, 4, 1, 4, \cdots \in \Sigma^+(X)$.  Choose $a_0 =\lfloor s - c\rfloor$.  Now $s - a_0 \in \eJ$.  Thus, by (\ref{eqn:HallThm}), there exists $\alpha:= [0; a_1, a_2, \cdots], \beta:= [0; b_1, b_2, \cdots]\in C_4$ such that $s-a_0 = \alpha+\beta$, giving $f(\theta^{(s)})= s$.  We claim that $f(\sigma^k \theta^{(s)}) \leq s$ for all $k$, which implies the desired result.  We now prove this claim.  Consider first $a_0 \geq 6$.  This condition gives $s > 6$.  Now, since $(\sigma^k \theta^{(s)})_0 \in X$ for $k \neq 0$, the claim follows.  
 
 Finally, consider $a_0=5$.  First suppose $a_1, b_1 \geq 2$.  For $k >0$, we have $f(\sigma^k \theta^{(s)}) = \theta^{(s)}_k + u + v$ where \[v:= [0;\theta^{(s)}_{k-1}, \theta^{(s)}_{k-2}, \cdots,  \theta^{(s)}_{1} , 5,  \theta^{(s)}_{-1},  \theta^{(s)}_{-2}, \cdots]\] and $u \in C_4$.  Set $\nu:=\theta^{(s)}_{k-1} \theta^{(s)}_{k-2} \cdots  \theta^{(s)}_{1}  5  \theta^{(s)}_{-1}  \theta^{(s)}_{-2} \cdots \in \Sigma^+(\{1, 2, 3, 4, 5\})$.  Let $\ell \in \NN$ be the least that $\nu_\ell \neq \chi_\ell$.  Consider first $\ell < k$.  If $\ell$ is even, then $\chi_\ell = 4 > \nu_\ell$, and, if $\ell$ is odd, then $\chi_\ell = 1 < \nu_\ell$.   In either case~\cite[Chapter~1, Lemma~1]{CF89} gives $r > v$.
 
 Now consider $\ell = k$.  If $k$ is odd, then $\chi_k = 1 < \nu_k=5$, giving $r > v$.  If $k$ is even, then $\nu_{k-1} =\theta^{(s)}_{1} = a_1\geq 2$, but $\chi_{k-1}=1$, which contradicts the choice of $\ell$.  Thus, this case cannot occur.  This shows $r > v$ and thus $f(\sigma^k \theta^{(s)}) < s$ for $a_1, b_1 \geq 2$, and $k >0$.  The analog proof shows $f(\sigma^k \theta^{(s)}) < s$ for $a_1, b_1 \geq 2$, and $k <0$.
 
 Next, suppose $a_1=1$.  Then $\alpha \geq (1 + r)^{-1}$ and thus $s \geq 5 + \min C_4 + (1 + r)^{-1} = (65 +11\sqrt 2)/14 > 4 + 2\max C_5$.  Now, $f(\sigma^k \theta^{(s)}) \leq 4 + 2\max C_5$ for every $k \neq 0$, giving the claim.  For $b_1=1$, the analogous proof gives the claim.
 
 %that $f(\sigma^k \theta^{(s)}) < 6$ for all $k \neq 0$.
\end{proof}

Now recall that $X_* := \{11, 22\}$ and every element of $H_Q(X_*)$ has a representative in $X_*^*$.  By Proposition~\ref{prop:XmTransitive}, $H_Q(X_*)$ acts transitively on $\PP^1(R_Q)$.  Thus, for every $\xi, \eta \in \PP^1(R_Q)$, there exists $w_{(\xi, \eta)} \in  X_*^*$ such that $\rho_Q(w_{(\xi, \eta)}) \xi =\eta$.  Note that $|w_{(\xi, \eta)}|_{\NN} < \infty$.  

%Define \[L_Q := \max_{\xi, \eta \in \PP^1(R_Q)} |w_{(\xi, \eta)}|_{\NN} < \infty.\] 

\begin{lemm}\label{lemm:HallSteering} Fix $Q \geq 2$, $s \geq s_\#$, and $[c_0;c_1, \cdots, c_\ell]$ where $c_0 \in \ZZ$, $\ell \in \NN \cup \{0\}$, and, when $\ell \neq 0$, $c_i \in \NN$ for all $1 \leq i  \leq \ell$.  Then there exists a sequence $(d_i)_{i = \ell+1}^\infty$ of positive integers such that $y := [c_0;c_1, \cdots, c_\ell, d_{\ell+1}, d_{\ell+2}, \cdots]$ with $k(y) = s$, and, for every $j \in R_Q$, there exists a strictly increasing sequence $(n_{r, j})_{r=1}^\infty$ in $\NN$ such that\begin{align}\label{eqn:Disp:lemm:HallSteering}
  [p_{n_{r, j}}: q_{n_{r, j}}]= [-j:j+1] \textrm { in } \PP^1(R_Q), \quad \lim_{r \rightarrow \infty} \Lambda_{n_{r, j}}(y) = s. \end{align}
 
\end{lemm}
\begin{proof}
 If $\ell = 0$, define $c_1 \cdots c_\ell$ to be the empty word.  Also, whenever $\epsilon_{\widetilde{r}}$, for some $\widetilde{r} \in \NN$, is used in the proof of this lemma, it comes from the analog of the proof of Lemma~\ref{lemm:transSWthres}.  Let $a_0$, $(a_i)$, $(b_i)$, and $\theta^{(s)}$ be as in Lemma~\ref{lemm:AttainHeightWithBigCenter}.  We adapt our distinguished word construction.  For all $r \in \NN$, define \[W^-_r := b_r \cdots b_1, \quad W^+_r := a_1 \cdots a_r, \quad W_r:= W^-_r a_0 W^+_r, \quad \xi_r:= \rho_Q(W^-_r) [0:1]\] where $a_0$ is referred to as the {\em distinguished digit (of $W_r$)}.  For $j \in R_Q$, let $\nu_j := [-j:j+1]$. 
 
 Let the strict ordering $\prec$ on $(r, j)$, where $r \in \NN$ and $j = 0, \cdots Q-1$, be the lexicographical, and let $\next(\cdot)$ be the successor function.  Let \[g_* := \begin{pmatrix}  1 & c_0\\ 0 & 1 \end{pmatrix} \in \PGL(2, R_Q)\] and define $g_{1,0}:= g_* \rho_Q(c_1
 \cdots c_\ell)$.  For $(r,j)$, define \[b_{r,j}:= w_{(\xi_r, g_{r,j}^{-1} \nu_j )}, \quad \B_{r,j}:= b_{r,j} W_r, \quad g_{\next(r,j)}:= g_{r,j} \rho_Q(\B_{r,j}).\]    
 
 The recursion yields $\gamma:= c_1 \cdots c_\ell \B_{1, 0} \cdots \B_{1, Q-1} \B_{2, 0} \cdots \in \Sigma^+(\NN)$.  Let $d_{\ell+i} = \gamma_{\ell+i} $ for all $i \in \NN$. Thus, we have $y = [c_0; c_1, \cdots, c_\ell, d_{\ell+1}, d_{\ell+2} \cdots]$.  Define \[n_{r, j}:=  \ell  +r +  |b_{r,j}|_{\NN} + \sum_{(r', j') \prec(r,j)} |\B_{r',j'}|_{\NN}. \]  Note that $[p_{n_{r, j}}: q_{n_{r, j}}]= g_{r,j} \rho_Q(b_{r,j}) \xi_r = \nu_j$ and that the distinguished digit $a_0$ associated with $(r,j)$ is at position $n_{r, j}+1$.  Thus, by Lemma~\ref{lemm:AttainHeightWithBigCenter}, we have $|\Lambda_{n_{r, j}}(y) - s| \leq \epsilon_r$.
 
 It remains to show that $k(y) \leq s$.  Fix $r_0 \geq 1$.  First, let \[\widetilde{N}_{r_0}:= \max_{\substack{1\leq r <2 r_0 \\ 0\leq j < Q}}n_{r,j}, \quad N_{r_0} := 1 + r_0  +\max (\ell,\widetilde{N}_{r_0}).\]  We claim that, for all $n \geq N_{r_0}$, if $|n -n_{r, j}| \leq r_0$, then $r \geq 2 r_0$.  We now prove the claim.  If $r < 2 r_0$, then $n \leq n_{r,j}+r_0 < N_{r_0}$, a contradiction.  
 
 There are two cases.  If positions $n+1 - r_0, \cdots, n+1+r_0$ do not contain any distinguished digit of the $W_r$ in any $\B_{r, j}$, then we have $\gamma_{n+1 - r_0}, \cdots, \gamma_{ n+1+r_0} \in X$.  Thus extending $\mu_{n+1 - r_0}:=\gamma_{n+1 - r_0}, \cdots, \mu_{ n+1+r_0}:=\gamma_{ n+1+r_0}$ to any element $\mu \in \Sigma(X)$ yields $f(\sigma^{n+1}\mu) \leq s_\#$ by Lemma~\ref{lemm:HallRayBasics} and, consequently, $\Lambda_n(y) \leq s_\# +  \epsilon_{r_0}$.
 
 Otherwise, positions $n+1 - r_0, \cdots, n+1+r_0$ contain the distinguished digit of the $W_r$ of some $\B_{r, j}$.  Hence, $r \geq 2r_0$ and  $|n -n_{r, j}| \leq r_0$. Since $|n -n_{r, j}|+r_0 \leq 2 r_0$, we have that the $\gamma_{n+1 - r_0} \cdots \gamma_{ n+1+r_0}$ is a subword of $W_r$.  Thus,  by Lemma~\ref{lemm:AttainHeightWithBigCenter}, we have \[\Lambda_n(y) \leq f(\sigma^{k}\theta^{(s)})+ \epsilon_{r_0}\leq s + \epsilon_{r_0}\] where $k:= n -n_{r, j}$.  Letting $r_0 \rightarrow \infty$ yields the desired result. \end{proof}

\begin{prop}\label{prop:HallRay}
Let $Q \geq 2$ be an integer and $F\subset \QQ$ be nonempty and finite with $Qr \in \ZZ$ for all $r \in F$.  Then $[4\sqrt 2 \, Q, \infty) \subset \eL(E_F)$.   Moreover, for every $t\geq4 \sqrt{2} \, Q$, the set \[\{x\in E_{\frac{1}{Q}\ZZ}:k(x)=t\}\]
is dense in $\RR$.
\end{prop}

\begin{proof} Fix $t \geq 4\sqrt 2 \, Q$ and put $s:= t/Q$.  Let $y$ be as in Lemma~\ref{lemm:HallSteering}.  The lemma gives that $k(y) = s \geq s_\#$ and (\ref{eqn:Disp:lemm:HallSteering}) holds for all $j \in R_Q$.  Now we may regard $QF \subset R_Q$.  Then Lemma~\ref{lemm:LagrangeUnderQPush} gives that $k(A_j(y)) = Qs=t$ for every  $j \in QF\cup \{0\}$.  Setting $x:= A_0(y)$ gives $t \in \eL(E_F)$ and thus yields the first assertion.

Let $I$ be a nonempty open interval in $\RR$. Choose $c_0, \cdots, c_\ell$ for Lemma~\ref{lemm:HallSteering} so that \[A_0([c_0; c_1, \cdots, c_\ell, c_{\ell+1}, c_{\ell+2}, \cdots]) \in I\] for any choice of sequence $(c_j)_{j= \ell+1}^\infty$ of positive integers.  Now, for these $c_0, \cdots, c_\ell$, Lemma~\ref{lemm:HallSteering}  constructs a $y=:y_I$ for which Lemma~\ref{lemm:LagrangeUnderQPush} gives that $k(A_j(y_I)) = Qs=t$ for every  $j \in R_Q$ (and thus every $j \in\ZZ$).  Note that $x_I = A_0(y_I) \in I$, which yields the second assertion.

\end{proof}

\section{The dimension of $E_F$ and low Lagrange equality spectra}\label{sec:DimEF} In this section, we prove our main result  Theorem~\ref{thm:EqualSetsBad} for the sets $E_F$ themselves.  To achieve this, we introduce fixed-target projective cocycle steering.  Arithmetic enters in the form of Proposition~\ref{prop:XmTransitive} and the following variant of Lemma~\ref{lemm:ExactCharMaxTransport}.

\begin{lemm}\label{lemm:GroundState}  For a reduced rational number $p/q \in \QQ$ such that $p+q \neq 0$, we have that \begin{align}\label{eqn1:lemm:GroundState}
[p:q]=[-1:1] \textrm { in } \PP^1(R_Q)  \implies e_j(p,q) =1 \textrm{ for all } j \in \ZZ. \end{align}  Moreover, let $x \in \RR \backslash \QQ$.  Then along any sequence of reduced rational numbers $p/q \rightarrow x$ such that $p+q \neq 0$ and $[p:q]=[-1:1]$ hold, we have that \begin{align}\label{eqn2:lemm:GroundState}
 \frac{\Lambda\left(A_j(x); \red\left(A_j(p/q) \right)\right)}{\Lambda(x;p/q)} = \frac 1 Q \frac{q|1+x|}{|p+q|}= \frac 1 Q(1 +o(1)). \end{align}

\end{lemm}

\begin{proof}  If $[p:q]=[-1:1]$, then $p+q \equiv 0 \pmod Q$, $p \in R^*_Q,$ and $e_j(p,q) = \gcd(p + j(p+q), Q(p+q)) = \gcd(p,q) =1$, giving (\ref{eqn1:lemm:GroundState}).  Applying Lemma~\ref{lemm:ApproxCoeffIntegralMob} gives (\ref{eqn2:lemm:GroundState}).
 
\end{proof}

\begin{prop}\label{prop:EFThick}
For every integer $Q \geq 2$, $\dim_H E_{\frac 1 Q \ZZ}=1$, and $ E_{\frac 1 Q \ZZ}$ is thick. Consequently, for every nonempty, finite $F \subset \QQ$, $\dim_H E_{F}=1$, and $ E_{F}$ is thick.
\end{prop}
\begin{proof}

Let $M \geq 2$ be an integer, $Y:=\{1, \cdots, M\}$, $0< u < \dim_H K(Y)$, and $H:= H_Q(Y)$.  Then, for $\varepsilon:= (\dim_H K(Y)- u)/2$, Lemma~\ref{lemm:IdentityAlphabet} gives a finite, prefix-free alphabet $A \subset Y^+$ with $L \in \NN$ such that, for every $a \in A$, we have \begin{align} |a|_{\NN}= L, \quad \rho_Q(a) = 1_H, \quad \dim_H K(A) >u.
  \end{align}

Now recall that $X_* := \{11, 22\}$.  By Proposition~\ref{prop:XmTransitive}, \begin{align}\label{eqn:MarkerFullHDEF} \textrm{ there exists } c_- \in X_*^+ \textrm{ such that }
 \rho_Q(c_-)[0:1] = [-1:1]. \end{align}  Let $b = Q^2(M+2) +1$.  Then $\rho_Q(b) = \rho_Q(1) \in H_Q(\{1, 2\})$.  Choose $c_+ \in \{1, 2\}^*$ such that $\rho_Q(c_+) = \left(\rho_Q(c_-) \rho_Q(b)\right)^{-1}$ and define \[W:= c_- b c_+\] and refer to $b$ as the {\em distinguished digit}.  Note that $b \notin \{1,2\}$ and $\rho_Q(W)=1_H$.

By Lemma~\ref{lemm:pressure}~\ref{it1:lemm:pressure}, for sufficiently large $m$, we have $Z^-_m(A, u) \geq e^{\eta m}$.  Define \[\Am := \{wW: w \in A^m\}, \quad \ell := m L + |W|_{\NN}\] for $m$ large enough so that Lemma~\ref{lemm:pressure}~\ref{it2:lemm:pressure} gives $\dim_H K(\Am) >u$.  This is analogous to the proof in Lemma~\ref{lemm:ConstXhatFromX}.

Let \[x:= [0; a_1, a_2, \cdots] \in K(\Am), \quad a_1 a_2 \cdots = w_1W w_2 W \cdots, \quad w_i \in A^m\] and $p_n/q_n$ be the convergents of $x$.  The digits equal to $b$ occur in positions $n_k +1$ where $n_k := k \ell + mL + |c_-|_{\NN}$ for $k \in \NN \cup \{0\}$, and hence $\rho_Q(a_1 \cdots a_{n_k}) = \rho_Q(c_-)$.  Consequently, we have \begin{align}\label{eqn:GroundStateCond}
 [p_{n_k}: q_{n_k}] = [-1:1], \end{align} which, by Lemma~\ref{lemm:GroundState}, gives \begin{align}\label{eqn:hitsFullDimEF}
 e_j(p_{n_k}, q_{n_k})=1 \end{align} for all $j \in \ZZ$ and all $k$.

Since all digits other than $b$ are in $Y$, (\ref{eqnPerron}) gives that \begin{align}\label{eqnMissConvFullDimEF}
 \Lambda_n(x) < M + 2  \end{align}for every $n \notin \{n_k: k \geq 0\}$.  At $n \in \{n_k: k \geq 0\}$, (\ref{eqnPerron}) gives \[b + 2/(1+b) < \Lambda_n(x)= b + [0; a_{n+2}, a_{n+3} \cdots] + [0; a_n, a_{n-1}, \cdots, a_1]< b+2 \] where the first inequality follows because $b >M$.  Let $\lambda := k(x)$.  Consequently, (\ref{eqnPerronConv}) gives that $\lambda = \limsup_{k \rightarrow \infty} \Lambda_{n_k}(x)$ and thus $\lambda/Q  > b/ Q \geq Q(M+2)>2Q$.

We claim that $k(A_j(x)) = \lambda/Q$ for all $j \in \ZZ$.  We now prove the claim. Choose $n_{k_i}$ such that $\Lambda_{n_{k_i}} (x)\rightarrow \lambda$, and let \[\frac{\widetilde{p}_i}{\widetilde{q}_i} := \red \left(A_j \left(\frac{p_{n_{k_i}}}{q_{n_{k_i}}} \right) \right), \quad \widetilde{q}_i >0.\]  Lemma~\ref{lemm:GroundState} and (\ref{eqn:GroundStateCond}) imply that $\widetilde{q}_i \rightarrow \infty$ and that $\Lambda(A_j x, \widetilde{p}_i/\widetilde{q}_i) \rightarrow \lambda/Q$.  Thus, \[k(A_j(x)) = \limsup_{q \rightarrow \infty} \frac 1 {q \|q A_j(x)\|} \geq \lim_{i \rightarrow \infty} \frac 1 {\widetilde{q}^2_i |A_j(x) - \widetilde{p}_i/\widetilde{q}_i|} = \frac \lambda Q,\] giving the lower bound.   Note that the $\widetilde{p}_i/\widetilde{q}_i$ are not necessarily convergents of $A_j(x)$.

We now prove the upper bound.  Recall the sets $\M_j(x)$ and $\U_j(x)$ of matched and unmatched convergents $P_\nu/Q_\nu$ from Section~\ref{subsec:MatchedUnmatchedConv}.  Also recall that, in taking the limit superior, a finite set of indices does not alter the limit superior and may thus be discarded.  Let \[\M_j^{\hit}(x) := \left\{ \nu \in \M_j(x) : P_\nu/Q_\nu = \red \left(A_j \left(\frac{p_{n_{k}}}{q_{n_{k}}}, \right) \right),  k \geq 0\right\}, \quad \M_j^{\miss}(x):=\M_j(x) \backslash \M_j^{\hit}(x).\]  Now for $\M_j^{\hit}(x)$, (\ref{eqn:GroundStateCond}) applies and thus Lemma~\ref{lemm:GroundState} gives \[\limsup_{\substack{\nu \rightarrow \infty \\ \nu \in\M_j^{\hit}(x) }} \Lambda_\nu (A_j(x)) \leq \frac \lambda Q,\] and, for $\M_j^{\miss}(x)$, Lemma~\ref{lemm:ApproxCoeffIntegralMob} and (\ref{eqnMissConvFullDimEF}) give \[\limsup_{\substack{\nu \rightarrow \infty \\ \nu \in\M_j^{\miss}(x) }} \Lambda_\nu (A_j(x)) \leq  Q(M+2).\]  Consequently,  Lemma~\ref{lemm:UnmatchedConvergents} implies that $k(A_j(x)) \leq \max(\lambda/Q, Q(M+2), 2Q) = \lambda/Q$, giving the claim.

The claim implies that $A_0(K(\Am)) \subset E_{\frac 1 Q \ZZ}$.  Since $A_0$ is bilipschitz, we have that $\dim_H E_{\frac 1 Q \ZZ}) >u$.  Now~\cite{Hen92} gives that $\dim_H K(Y) \rightarrow 1$ as $M \rightarrow \infty$.  By letting $u \rightarrow \dim_H K(Y)$ and then letting $M \rightarrow \infty$, we have $\dim_H E_{\frac 1 Q \ZZ} =1$.

Let $I$ be a nonempty open interval in $\RR$.  For $d_0 \in \ZZ$ and $d_1, \cdots, d_t \in \NN$ define \[\mP := \begin{pmatrix}  1 & d_0 \\ 0 & 1 \end{pmatrix} M_{d_1 \cdots d_t}.\]  For the thickness assertion, choose $d_0 \in \ZZ$ and $d_1, \cdots, d_t \in \NN$ so that $A_0(\mP([0,1])) \subset I$. Now replace (\ref{eqn:MarkerFullHDEF}) with  \begin{align*}\label{eqn:MarkerthickEF} \textrm{ there exists } c_- \in X_*^+ \textrm{ such that }
 \mP\rho_Q(c_-)[0:1] = [-1:1]. \end{align*}  With this replacement, the analog of the above proof gives $ \dim_H \left(I \cap E_{\frac 1 Q \ZZ}\right) =1$, showing the thickness assertion.

\end{proof}

\begin{prop}\label{prop:EQQIsEmpty}
Let $x \in \RR \backslash \QQ$.  Then $\sup_{r \in \QQ} k(x+r) = \infty$.
\end{prop}

\begin{proof}
 It suffices to show \[\max_{0 \leq j <p} k(x + \frac j p) \geq p\] for every prime $p$.  Fix a prime $p$ and let $P_n/Q_n$ be the convergents of $px$.  Define $\mN_p := \{n : p \nmid Q_n\}$.  Since $\gcd(Q_n, Q_{n+1})=1$, $|\mN_p|=\infty.$  Now for each $n \in \mN_p$ there exists a unique $j_n \in R_p$ such that $P_n + j_n Q_n \equiv 0 \pmod p$.  By the pigeonhole principle, there exist $j \in R_p$ and infinite subset $\mN \subset \mN_p$ such that $j_n = j$ for all $n \in \mN$.
 
 Now define the integer $a_n := (P_n + j Q_n)/p$.  Let $d \in \NN$ be such that $d \mid a_n$ and $ d \mid Q_n$.  Then $d \mid (p a_n - j Q_n)$.  Since $P_n = p a_n - j Q_n$ and $\gcd(P_n, Q_n)=1$, we have that $d = 1$.
 
 Thus, $Q_n^2 | x + j/p - a_n/Q_n| = Q_n^2/p |px - P_n/Q_n| < 1/p$ for all $n \in \mN$.  Taking the limit superior as $Q_n \rightarrow \infty$ along $n \in \mN$ gives that $k(x+j/p) \geq p$. 
\end{proof}

\begin{coro}\label{cor:EQQIsEmpty}
We have that $E_\QQ = \emptyset$. 
\end{coro}

\begin{proof}
 Assume not.  Let $x \in E_\QQ$.  Then we have $\sup_{r \in \QQ} k(x+r) = k(x) < \infty$ because $x \in \BA$, contradicting Proposition~\ref{prop:EQQIsEmpty}.
\end{proof}

\begin{coro}\label{cor:LowSpectra}
Let $Q \geq 2$ be an integer and let $p_Q$ be the largest prime dividing $Q$.  Then \[\eL(E_{\frac 1 Q \ZZ}) \cap (-\infty, p_Q) = \emptyset.\]
\end{coro}
\begin{proof}
 Immediate from Proposition~\ref{prop:EQQIsEmpty}.
\end{proof}

\begin{coro}\label{cor:QisPStatement}  Let $p$ be prime.  We have \[[4\sqrt 2 p, \infty) \subset  \eL\left(E_{\frac 1 p \ZZ}\right) \subset [p, \infty), \quad p \leq \inf  \eL\left(E_{\frac 1 p \ZZ}\right) \leq 3p.\]
\end{coro}
\begin{proof}
 Apply Propositions~\ref{prop:HallRay} and~\ref{prop:PosHD} and Corollary~\ref{cor:LowSpectra}.
\end{proof}

\section{Simultaneous rational-translation Lagrange equality spectra for $\BA(\delta)$}  In this section, we state and prove the analogs of Theorem~\ref{thm:MainBad} and ~\ref{thm:EqualSetsBad} for $\BA(\delta)$.

\begin{theo}\label{thm:MainBadDelta}  Let $\delta>0$, $Q \geq 2$ be an integer,  $F:= \frac 1 Q \ZZ$, and $p_Q$ be the largest prime dividing $Q$.  Then we have \begin{enumerate}[label=(\roman*)]
\item\label{it1:HallRay:thm:MainBadDelta} $[4\sqrt 2 \, Q, \delta^{-1}) \subset \eL(E_F(\delta))$,
\item\label{it2:FullHD:thm:MainBadDelta}  $d_{F, \delta}(t) = 1$ for all $\delta^{-1} > t \geq Q \Ft_1$,
\item\label{it3:Scaling:thm:MainBadDelta}  $d_{F, \delta}(t) \geq d_{\cl}(t/Q)$ for all $\delta^{-1} > t \geq Qs_*$,
\item\label{it4:PosDim:thm:MainBadDelta}  $d_{F, \delta}(t) >0$ for all $\delta^{-1}> t > 3Q$,
\item\label{it5:LowSpec:thm:MainBadDelta}  $\eL(E_F(\delta)) \cap (-\infty, p_Q) = \emptyset$.

\end{enumerate}
Moreover, \begin{enumerate}[label=(\alph*)]
\item\label{it6:DenE:thm:MainBadDelta}  for every $\delta^{-1} > t\geq4 \sqrt{2} \, Q$, we have $\overline{\{x\in E_F(\delta):k(x)=t\}}= \RR$,
\item\label{it7:FullDimCpctSet:MainBadDelta} if $\delta^{-1} > 2 Q \sqrt 3$, then, for every $\varepsilon >0$, there exists a compact set $V_{Q, \varepsilon}$ with $\dim_H V_{Q, \varepsilon}=1$ such that \[V_{Q, \varepsilon} \subset \eL(E_F(\delta)) \cap (2 Q \sqrt 3 - \varepsilon, 2 Q \sqrt 3),\]
\item for a prime $p < (4\sqrt 2\delta)^{-1}$, we have \[[4\sqrt 2 p, \delta^{-1} ) \subset  \eL\left(E_{\frac 1 p \ZZ}(\delta) \right) \subset [p, \delta^{-1}], \quad p \leq \inf  \eL\left(E_{\frac 1 p \ZZ}(\delta)\right) \leq 3p.\]

\end{enumerate}
\end{theo}
\begin{proof}
% Apply Proposition~\ref{prop:BadDeltaRestriction}
 
 Apply (\ref{eqn:BadDeltaEndPointBnds}) to the corresponding parts of Theorem~\ref{thm:MainBad}.
\end{proof}

Let $M \in \NN$.  Recall $C_M$ from Section~\ref{sec:HallRay}.

\begin{theo}\label{thm:EqualSetsBadDelta}  Let $\delta>0$. \begin{enumerate}[label=(\roman*)]
\item\label{it1:thickEFinite:thm:EqualSetsBadDelta} Let $Q \geq 2$ be an integer, $M:= \lceil 1/Q\delta - 3/Q^2\rceil-3$, $N:=\lceil 1/\delta \rceil - 1$ and $I \subset \RR$ be an nonempty open interval.  If $\delta < Q/(4Q^2 +3)$, then \[\dim_H C_M \leq \dim_H \left(E_{\frac 1 Q \ZZ}(\delta) \cap I \right)\leq \dim_H C_N.\]
\item\label{it2:emptyEQ:thm:EqualSetsBadDelta} $E_\QQ(\delta) = \emptyset$.
\end{enumerate}
\end{theo}
\begin{proof} The proof of the assertion~\ref{it1:thickEFinite:thm:EqualSetsBadDelta} is as follows.  First note that $M\geq 2$.  Fix $0< u  < \dim_H C_M$ and $b= Q^2(M+2) +1$.  Then the proof of Proposition~\ref{prop:EFThick} gives a compact set $K \subset I$ such that, for all $x \in K$ and $j \in \ZZ$, we have that $k(x+j/Q) = k(x) \leq (b+2)/Q < \delta^{-1}$.  Using (\ref{eqn:BadDeltaEndPointBnds}) and letting $u \rightarrow \dim_H C_M$ give the lower bound.  For the upper bound, note that the elements of $\BA(\delta)$ have partial quotients that are eventually bounded by $N$, giving $\dim_H \BA(\delta) \leq \dim_H C_N$.  Finally, the assertion~\ref{it2:emptyEQ:thm:EqualSetsBadDelta} follows from Theorem~\ref{thm:EqualSetsBad}~\ref{it2:emptyEQ:thm:EqualSetsBad} because $E_\QQ(\delta) \subset E_\QQ$ by (\ref{eqn:BadDeltaEndPointBnds}).
\end{proof}

\end{document}